\documentclass[12pt]{svj3}                     
\smartqed  
\usepackage{graphicx}

\usepackage{amsmath}
\usepackage{amsfonts,amssymb}
\newcommand{\doc}{\underline{\it Proof.}\
}
\newcommand{\bo}{\qed}
\newcommand{\rav}{\stackrel{\triangle}{=}}
\newcommand{\ph}{\varphi}
\newcommand{\la}{\langle}
\newcommand{\ra}{\rangle}
\newcommand{\epsi}{\varepsilon}
\newcommand{\beg}{\varkappa}

\newcommand{\rref}[1]{$(\ref{#1})$}
\newcommand{\mm}[1]{{\bf{#1}}}
\newcommand{\ct}[1]{{\mathcal{#1}}}
\newcommand{\fr}[1]{{\mathfrak{#1}}}
\newcommand{\td}[1]{\widetilde{#1}}
\newcommand{\ravref}[1]{\stackrel{(\ref{#1})}{=}}
\newcommand{\leqref}[1]{\stackrel{(\ref{#1})}{\leq}}
\newcommand{\geqref}[1]{\stackrel{(\ref{#1})}{\geq}}
\newcommand{\sref}[1]{\stackrel{(\ref{#1})}{\vartriangleright}}
\usepackage[active]{srcltx}
\newcommand{\ee}{
}

\begin{document}

\title{THE TRANSVERSALITY CONDITIONS IN  INFINITE HORIZON PROBLEMS
AND THE STABILITY OF ADJOINT VARIABLE\thanks{Our grant could be
here}}

\titlerunning{NECESSARY CONDITIONS IN CASE OF STABLE ADJOINT VARIABLE}     

\author{Dmitry Khlopin 
}

\institute{D. Khlopin \at
              Institute of Mathematics and Mechanics, Ural Branch,   Russian
Academy of Sciences, Yekaterinburg, Russia \\
              \email{khlopin@imm.uran.ru}           
}

\date{Received: 10.08.2011 / Accepted: date}

\maketitle

\begin{abstract} This paper investigates the
necessary conditions  of optimality for uniformly overtaking optimal
control  on infinite horizon with free right endpoint. Clarke's form
of the Pontryagin Maximum Principle is proved without the assumption
on boundedness of total variation of adjoint variable. The
transversality condition for adjoint variable is shown to become
necessary if the adjoint variable is partially Lyapunov stable. The
modifications of this condition are proposed for the case of
unbounded adjoint variable. The Cauchy-type formula for the adjoint
 variable proposed
  by S.\,M.~Aseev   and A.\,V.~Kryazhimskii
in~\cite{kr_as1},\cite{kr_as2} is shown to complement relations of
the Pontryagin Maximum Principle up to the complete set of necessary
conditions of optimality if the improper integral in the formula
converges conditionally and continuously depends on the original
position. The results are extended to an unbounded objective
functional (described by a nonconvergent improper integral),
unbounded constraint on the control, and uniformly sporadically
catching up optimal control.  \keywords{Optimal control \and
infinite horizon problem \and transversality condition for infinity
\and necessary conditions \and Lyapunov stability  \and uniformly
overtaking optimal}
 \subclass{ 49K15 \and 49J45 \and 37N40 \and 34D05}
\end{abstract}

\section*{Introduction}
\label{intro} The Pontryagin Maximum Principle for infinite
horizon problems has already been formulated in the
monograph~\cite{ppp}, but without the transversality condition the
obtained relations were incomplete and in general, selected a much
too broad family of potentially extremal trajectories.
A significant number \cite{Halkin,kr_as,aucl,mich,smirn,sch,kam} of
such conditions has been proposed; however, as it was noted in, for
example, \cite{Halkin,mich,shell},\cite[Sect. 6]{kr_as},\cite[Example
10.2]{norv}, these conditions may be either inconsistent with the
relations of the Pontryagin Maximum Principle, or follow from them.
Hence the need to investigate the applicability of a transversality
condition (see~\cite{kr_as,aucl,kam,mich,smirn,Ye,norv,ss,sagara})
and the need to separately check if it is necessary for a specific
optimization problem. The first aim of this paper is to offer a
common approach to selecting a necessary transversality condition on
the adjoint variable for this problem
(Subsect.~\ref{sec:34}--\ref{sec:35}). However, the necessity of a
condition does not imply its nontriviality on solutions of the
relations of the Maximum Principle. Hence the need to find a
condition that would select a single solution of the relations of the
Maximum Principle for any uniformly overtaking optimal control. In
the papers \cite{kr_as1,kr_as2,kr_asD,kr_as,kr_as3}, Aseev and
Kryazhimskii develop and investigate the Cauchy-type formula for the
adjoint variable that possesses such a property. The second aim is to
maximize the applicability of the approach~\cite{kr_as}
(Sect.~\ref{sec:6}).

First of all,
 we construct the bicompact extension (see~\cite{va})
for the space of admissible controls in the form of the inverse limit
of the sequence of corresponding finite horizon extensions.
 It is shown that there exists a uniformly overtaking optimal generalized control for the case of a conditionally convergent objective functional that converges uniformly  with respect to all trajectories; this generalizes some results \cite{bald},\cite{carm},\cite{dm}.
Without this assumption, for uniformly overtaking optimal control for problems with
free right endpoint, the necessity of the Pontryagin Maximum
Principle in Clarke's
 form for the more general
conditions than in \cite{kr_as,kr_as3,mich},\cite[Theorem 2.1]{Ye} is shown; the obtained result is not a part of the results~\cite{av,Halkin,sagara}.

In Subsect.~\ref{sec:23} for a free right endpoint problem, the
convergence of transversality conditions on the adjoint variable is
provided by the integral partial stability of the adjoint variable as
a component of the Maximum Principle system. Thus we obtain the
assumptions that guarantee the necessity of such condition, which are
more general than the assumptions in \cite{norv},\cite[Corollary
2.1]{Ye},\cite[Theorem 1]{av} (It seems that the first implementation
of the approach that employs the notion of stability to obtain
transversality conditions was in~\cite{smirn}).
 If the the objective function and the right-hand side of the
 equation of dynamics are smooth (in the phase variable),
 then, instead of integral partial stability we can check
  the simpler condition of partial Lyapunov stability for
  the variable~$\psi$ as a component of solutions of the
  system of the Maximum Principle. For example, we can check
  if all Lyapunov exponents are negative for this variable 
(see Subsect.~\ref{sec:34}).

We propose the new transversality condition: the product of the
adjoint variable and a matrix function of time must be vanishingly
small at infinity. This condition becomes necessary if the product is
stable. The stability can be provided by the correct choice of the
matrix function; the choice may also reflect a priori information on
stability and asymptotic estimates, which may allow to pinpoint the
single extremal (see Subsect.~\ref{sec:35}).

If the above matrix function is the fundamental matrix of linearized
system
along the optimal trajectory, then, the corresponding transversality
condition automatically yields the formula that was proposed for
affine systems in~\cite{kr_as1},\cite{kr_as2}, general case of which
was examined in~\cite[Theorems~11.1, 12.1]{kr_as}, \cite[Theorem
2]{kr_as3},\cite[Theorem 1]{av}. As it was shown
in~\cite[Sect.~16]{kr_as}, the results of \cite{aucl,Ye} are the
corollaries of \cite[Theorem 12.1]{kr_as},\cite[Theorem 2]{kr_as3}.

Such choice of matrix function allows us to reduce the question
of necessity of the corresponding condition to not just the question
of the stability of the product, but even to the issue of checking if
the improper integral from the formula~\cite[(12.8)]{kr_as},
\cite[Theorem 2]{kr_as3},\cite[Theorem 1]{av} converges conditionally
and continuously depends on the initial position of the original
problem.
 This yields the Cauchy-type formula for the
adjoint variable and the ``normal'' Pontryagin Maximum Principle
under the assumptions weaker than in~\cite[Theorem~12.1]{kr_as}. This
result also generalizes~\cite[Theorem
3.2]{sagara},\cite[Theorem~2.1]{Ye} and~\cite[Theorems~3.1
and~8.1]{norv} (as far as the necessary conditions for problems with
free right endpoint are concerned).

 For the case of monotonous system, we also demonstrated certain estimates for the adjoint variable. In particular, we obtained the nonnegativity of adjoint variable under the weaker assumptions than in~\cite[Theorem 1]{kr_as_t}, \cite[Theorem 1]{weber},
   \cite[Theorem 10.1]{kr_as}.

In the last part of the paper, we extend the obtained results to the the cases of $\sigma$-compact constraints on controls and
to uniformly sporadically catching up optimal controls. Important breakthroughs for these problems were recently achieved in~\cite{av}.

A part of the results of this paper has been shown and announced in paper~\cite{my1}.

\section{Preliminaries}
\label{sec:1} \ee
 We consider the time interval
  $\mm{T}\rav
  \mm{R}_{\geq 0}.$
The phase space of the original control system is the certain
finite-dimensional metric space $\mm{X}\rav\mm{R}^m$. The
 unit ball of this space is denoted by~$\mm{D}$. Let $\mm{L}$ denote the linear space of all $m\times
m$ matrices. For the sake of definiteness, let us equip~$\mm{L}$ with
the operator norm. The
symbol~$E$ (which may be equipped with some indices) denotes various
auxiliary finite-dimensional Euclidean spaces, and the symbol
$\ct{B}(E)$ denotes the $\sigma$-algebras of their Borel subsets.

For a subset~$A$ of a topological space, $cl\, A$ denotes the closure
of this subset.

On the sets of all functions that are continuous on the
whole~$\mm{T}$, we consider the topology of uniform convergence
on~$\mm{T}$ and the compact-open topology; for example,
$C(\mm{T},E)$ and $C_{loc}(\mm{T},E)$. The first one is considered
to be equipped with the norm  $||\cdot||_C$ of the uniform
convergence topology. $\Omega$~denotes the family of functions
$\omega\in C(\mm{T},\mm{T})$ such that $\displaystyle
\lim_{t\to\infty}\omega(t)=0.$

Here and below, for each summable function~$a$ of time, the integral
$\int_{\mm{T}}a(t)dt$ is the limit $\int_{[0,T]}a(t)dt$ as
$T\to\infty$. The integral over an infinite interval, for example,
over $[T,\infty\ra$, is interpreted in the same sense.

Let us also consider a finite-dimensional Euclidean space~$\mm{U}$
and a set-valued map $U:\mm{T}\rightsquigarrow\mm{U}.$ The
set~$\fr{U}$ of admissible controls is understood as the set of all
Borel measurable selectors of the multi-valued map~$U$. The topology
on~$\fr{U}$ is defined by virtue of the inclusion
  $\fr{U}\subset \ct{L}^1_{loc}(\mm{T},\mm{U}).$

A function $a:\mm{T}\times E'\times \mm{U}\to E''$ is said to
\begin{description}
\item[1)] satisfy the Carath\'{e}odory conditions if
 a) the function $a(\cdot,y,u):\mm{T}\to E''$ is measurable for all $(y,t,u)\in\mm{X}\times Gr {U},$
 b) the function $a(t,\cdot,\cdot):E'\times U(t)\to E''$ is continuous for all $t\in\mm{T}.$
\item[2)]  be locally Lipshitz continuous
  if for each compact $K\in(comp)(\mm{T}\times E)$ there exists
  a function $L_K^a\in \ct{L}^1_{loc}(\mm{T},\mm{T})$ such that for all
  $(t,x'),(t,x'')\in K,u\in U(t)$, the inequality
  $||a(t,x',u)-a(t,x'',u) ||_{E''}\leq L_K^a(t)||x'-x''||_{E'}$ holds.
\item[3)]  be 
 integrally bounded (on each compact subset of $\mm{T}\times E$)
 if for each compact $K\in(comp)(\mm{T}\times E)$ there exists a function
 $M_K^a\in \ct{L}^1_{loc}(\mm{T},\mm{T})$ such that for all
  $(t,x)\in K,u\in U(t)$ we have $||a(t,x,u)||_{E''}\leq M_K^a(t)$.
\item[4)] satisfy the continuability condition on~$\mm{T}$ if it satisfies
the sublinear growth condition, i.e., if the function~$f$
is Lipshitz continuous such that the
 function $L_K^{a}$ is independent of~$K$ and is
 integrally bounded (on each compact subset);
  see~\cite[1.4.6]{tovst1}.
\end{description}

\medskip

 Here and below, we assume the following conditions hold:

 {\it \bf Condition}~${\bf{(u)}}:$
    $U$ is a compact-valued map such that it is integrally
  bounded (on each compact subset of~$\mm{T}$) and 
  $Gr\,  U  \in       \ct{B}(\mm{T}\times\mm{U})$. 

 {\it \bf Condition}~${\bf{(fg)}}:$ the mappings $f:\mm{T}\times\mm{X}\times \mm{U}\to
 \mm{X}$,$g:\mm{T}\times\mm{X}\times \mm{U}\to\mm{R}$
 are locally Lipshitz continuous Carath\'{e}odory mappings that are integrally bounded
 (on each compact subset) and $f$ satisfies the continuability condition.

 Let us
consider the control system
\begin{subequations}
   \begin{equation}
   \label{sys}
    \dot{x}=f(t,x,u),\  x(0)=0,\qquad t\in\mm{T},\ x\in\mm{X},\ u\in U(t).
   \end{equation}
Now we can assign the solution \rref{sys} to every~$u\in\fr{U}$. The
solution is unique and it can be extended to the whole $\mm{T}$. Let
us denote it by~$\ph[u]$. 
The mapping $\ph:\fr{U}\to C_{loc}(\mm{T},\mm{X})$ is continuous.

In what follows, we examine the problem of maximizing the objective functional
\begin{eqnarray}
   \lim_{T\to\infty} J_T(u)\to\max;
    \nonumber \\[-1.5ex]
   \label{opt}
\\[-1.5ex]
   J_T(u)\rav\int_0^T g\big(t,\ph[u](t),u(t)\big)dt.
       \nonumber
   \end{eqnarray}
\end{subequations}
 If there is no  limit in~\rref{opt}, the optimality may be defined in diverse ways (for details, see~\cite{car},\cite{car1},\cite{stern}),
generally, we will us the following one:
\begin{definition}  A control $u^0\in\fr{U}$ is called
 uniformly  overtaking  optimal if for  each  $\epsi \in\mm{R}_{>0}$
    there
   exists  $T\in\mm{R}_{>0}$  such  that
 $
   J_t(u^0)\geq J_t(u)-\epsi $
 holds for all $u\in\fr{U}$, $t\in[T,\infty\ra$.
\end{definition}
Note that in paper~\cite{stern}, this definition is referred to as
uniformly  catching up optimal control. In \cite[Theorem~3.1]{2}, it
is shown that the uniformly overtaking optimality is equivalent to
the condition $$
   \lim_{t\to +\infty} \Big(J_t(u^0) - \sup_{u\in\fr{U}} J_t(u)\Big)=0$$
which, in terms of~\cite{2}, says that $u^0$ is
strongly agreeable.

   Note that for each uniformly  overtaking  optimal control $u^0\in\fr{U}$
 there exists a
function $\omega^0\in\Omega$ such that
\begin{equation}
 \label{555}
J_t(u^0)\geq J_t(u)-\omega^0(T)\qquad\forall u\in
 \fr{U},T\in\mm{T},t\in[T,\infty\ra.
   \end{equation}

\section{On existence of  uniformly  overtaking  optimal control}
\label{sec:11} \ee
To complete the main objective of this section, we need the following assumption:

 {\it \bf Condition}~${\bf{(e)}}:$
 there exists a
function $\omega\in\Omega$ such that
 $$\int_T^\tau g\big(t,\ph[u](t),u(t)\big) dt \leq \omega(T)\qquad\forall u\in
 \fr{U},T,\tau\in\mm{T}, T<\tau.$$

 Note that to the best of author's knowledge, a one-sided condition like ${\bf{(e)}}$ was first proposed in paper~\cite[$(\Pi 8)$]{dm}. As it was actually proved in~\cite[Subsect 5.1]{dm},
instead of~${\bf{(e)}}$, it is enough to assume, for example, the stronger condition
 $$g\big(t,\ph[u](t),u(t)\big)\leq l(t) \qquad\forall u\in
 \fr{U},T\in\mm{T}$$
for some summable  
 on  $\mm{T}$ mapping $l\in \ct{L}^1(\mm{T},\mm{R}).$

\subsection{The definition of  the set $\td{\fr{U}}$ of generalized
controls}

For each $u\in\mm{U}$, the symbol $\td\delta(u)$ denotes the
probability measure concentrated at the point $u$. Let
$\td{\fr{U}}_n$ denote the family of all weakly measurable
mappings~$\mu$ from $[0,n]$ to the set of Radon probability measures
over~$ \mm{U}$ such that $\int_{U(t)} \eta(t)(du)=1$ for a.a.
$t\in[0,n]$. Let us equip this set with the topology of *-weak
convergence. Then, the obtained topological space is a
compact~\cite[IV.3.11]{va}, and the set
$\fr{U}_n\rav\{u|_{[0,n]}\,|\,u\in\fr{U}\}$ is everywhere densely
included in $\td{\fr{U}}_n$ \cite[IV.3.10]{va} by the mapping
$u\to\td{\delta}\circ u$.

Now, let us introduce the set of all maps~$\eta$ from~$\mm{T}$ into
the set of Radon probability measures over~$ \mm{U}$ such that
$\eta|_{[0,n]}\in \td{\fr{U}}_n$ for every $n\in\mm{N}$; and let us
denote it by $\td{\fr{U}}$. To each~$n\in\mm{N}$ let projections
$\td\pi_n:\td{\fr{U}}\to \td{\fr{U}}_n$  be given by
$\td{\pi}_n(\eta)\rav\eta|_{[0,n]}$ for all $\eta\in\td{\fr{U}}$. Let
us equip~$\td{\fr{U}}$ with the weakest topology such that all
projections are continuous. The set~$\td{\fr{U}}$ is called the set
of generalized controls.

\medskip

Let us assume that for the certain Euclidean space~$E$ a mapping
$a:\mm{T}\times{E}\times \mm{U}\to
  (comp)(E)$ is given and the following condition is satisfied:

 {\it \bf Condition~}${\bf{(a)}}:$ the mapping
 $a:\mm{T}\times{E}\times \mm{U}\to (comp)(E)$ is a locally Lipshitz continuous
  integrally bounded Carath\'{e}odory mapping that satisfies the continuability condition.

Let us fix the set  $\Xi\subset E$ of initial values and the system
for $u\in\fr{U}$:
\begin{equation}
   \label{a}
   \dot{y}=a(t,y(t),u(t)),\  y(0)=\xi\in \Xi,\qquad
   t\in \mm{T}, u\in\fr{U}.
\end{equation}
It can also be generalized for $\eta\in\td{\fr{U}}$:
\begin{equation}
   \label{1650}
    \dot{y}=\int_{U(t)}a(t,y(t),u)\eta(t)(du),\ y(0)\in \Xi,\qquad
   t\in \mm{T}, \eta\in\td{\fr{U}}.
\end{equation}
Each its local solution can be extended to the whole~$\mm{T}$. For
every $\eta\in\td{\fr{U}}$, let us denote the family of all solutions
$y\in C_{loc}(\mm{T},E)$ of system~\rref{1650} by
$\td{\fr{A}}[\eta]$.

\subsection{The relaxed infinite-horizon optimal control problem.}

Similarly, we can consider the solution $\td\ph[\eta]\in
C_{loc}(\mm{T},\mm{X})$ of the Cauchy problem
\begin{subequations}
\begin{equation}
\dot{x}=\int_{U(t)} f(\tau,x(\tau),u)\,\eta(t)(du),\
x(0)=0\qquad \forall\eta\in \td{\fr{U}},
   \label{sysg}
   \end{equation}
the function
  $T\mapsto
  \td{J}_T(\eta)\rav \int_{[0,T]}\int_{U(t)} g(t,\td\ph[\eta](t),u)\,\eta(t)(du)\,dt;$
and
  the problem of maximizing the functional
\begin{equation}
   \lim_{T\to\infty} \td{J}_T(u)\to\max.
   \label{optg}
   \end{equation}
\end{subequations}
\begin{proposition}
 \label{spectr}
Assume ${\bf{(u)}}$. Then,\
\begin{description}
\item[1)] the space~$\td{\fr{U}}$ is a compact, and $\td\delta({\fr{U}})$ is everywhere dense in it;
\item[2)] If ${\bf{(a)}}$ holds, then for a compact $\Xi\in(comp)(E)$
   the map $\td{\fr{A}}:\td{\fr{U}}\to C_{loc}(\mm{T},E)$ is continuous, and
     $\td{\fr{A}}[\td{\delta}\circ\fr{U}]$ is everywhere dense in $\td{\fr{A}}[\td{\fr{U}}]\in
     (comp)(C_{loc}(\mm{T},E));$
\item[3)] If ${\bf{(fg)}}$ hold, then
  $\td\ph,\td{J}\in C_{loc}(\mm{T}\times\td{\fr{U}},\mm{R})$;
 \item[4)] If  ${\bf{(fg),(e)}}$ hold, then \
   there is a uniformly  overtaking  optimal  control $\td{u}^0\in\td{\fr{U}}$
   for the relaxed problem \rref{sysg}--\rref{optg} such that
\begin{eqnarray}
   \lim_{T\to\infty} \sup_{u\in {\fr{U}}} \int_0^T g(t,{\ph}[u](t),u(t))dt
   &=&
    \lim_{T\to\infty}\max_{\eta\in\td{\fr{U}}}\td{J}_T(\eta)=
    \max_{\eta\in\td{\fr{U}}}\lim_{T\to\infty} \td{J}_T(\eta)=
    \nonumber \\[-1.5ex]
   \label{supmax}
\\[-1.5ex]
    =
    \lim_{T\to\infty}\td{J}_T(\td{u}^0)&=&
 \int_{\mm{T}}\int_{U(t)} g(t,\td{\ph}[\td{u}^0](t),u)\td{u}^0(du)dt,    \nonumber
   \end{eqnarray}
and all limits in \rref{supmax} exist, although they can equal $-\infty$
\end{description}

\end{proposition}
\doc
For the sake of brevity, let us denote $\td\Pi\rav\prod_{n\in\mm{N}}
\td{\fr{U}}_n$, and let us equip it with Tikhonov topology.
 Let  $\td{\Delta}:\td{\fr{U}}\to\td\Pi$  be given by
 $\td{\Delta}(\eta)\rav\big(\td{\pi}_n(\eta)\big)_{n\in\mm{N}}$
for all $\eta\in\td{\fr{U}}.$ It is a homeomorphism by continuity of
the maps~$\td{\pi}_{n}$ and $\td{\pi}_n\circ\td{\Delta}^{-1}$.

Let $n,k\in\mm{N}, (n>k).$ Then, the space $\td{\fr{U}}_n$ is
included in~$\td{\fr{U}}_k$ by the mapping
  $\td{\pi}^{n}_{k}(\eta)\rav\eta|_{[0,k]}$ for all $\eta\in \td{\fr{U}}_n.$
 By $\td{\pi}^{n}_{k}\circ\td{\pi}^{k}_{i}=\td{\pi}^{n}_{i}$ for all
  $n,k,i\in\mm{N}, (n>k>i),$ we have the projective sequence of the topological spaces
  $
  \{
  \td{\fr{U}}_n,\td{\pi}^{n}_{k}\};$
and we can define the inverse limit \cite[III.1.5]{phillvv},
  \cite[2.5.1]{en}. In our notation, we can write it in the form
    $\lim_{\leftarrow} \{\td{\fr{U}}_n,\td{\pi}^{n}_{k}\} \rav\td\Delta(\td{\fr{U}})\subset \td\Pi.$
As shown above,~$\td\Delta$ is a homeomorphism; hence, $\td{\fr{U}}$
is homeomorphous to $\td\Delta(\td{\fr{U}})$. Now, by Kurosh
Theorem~\cite[III.1.13]{phillvv}, the inverse limit
$\td\Delta(\td{\fr{U}})$ of compacts $\td{\fr{U}}_n$ is compact, and
 $\td{\fr{U}}$ is a compact too. Similarly, from \cite[4.2.5]{en} and \cite[IV.3.11]{va} it
follows that $\td{\fr{U}}$ is also metrizable.

Repeating the reasonings without $\td{\ }$ or referring
to~\cite[3.4.11]{en} and~\cite[2.5.6]{en} yields
$\displaystyle{\fr{U}}\cong\lim_{\leftarrow} \{
  {\fr{U}}_n,\pi^{n}_{k}\}\rav \Delta({\fr{U}})\subset \Pi.$

For each $n\in\mm{N}$, let the mapping ${e}_n:\fr{U}_n\to
\td{\fr{U}}_n$ be given by
    ${e}_n(u)(t)\rav(\td\delta\circ u)(t)=\td{\delta}_{u(t)}$ for all
    $t\in[0,n],u\in\fr{U}_n.$ Since for all $n,k\in\mm{N},
  n>k$ it holds that
  ${e}_k \circ {\pi}^n_k = {e}_n,$
we have the projective system $\{{e}_n,{\pi}^n_k\}$. Passing to the
inverse limit, we obtain the mapping
$e_\Delta:\Delta(\fr{U})\to\td\Delta(\td{\fr{U}})$; from
$e_n\circ\pi_n=\td\pi_n\circ\td\delta$ we have
    $e_\Delta\circ\Delta=\td\Delta\circ\td{\delta}$, and from
  $\td{\fr{U}}_n=cl e_n(\fr{U}_n)$ (\cite{va}) we have
  $\td\Delta(\td{\fr{U}})=cl e_\Delta\big(\Delta(\fr{U})\big)=cl (\td\Delta\circ\td{\delta})(\fr{U});$
now, by continuity of~$\td\Delta^{-1}$, we obtain $\td{\fr{U}}= cl\,
\td{\delta}(\fr{U}).$

The mapping $\td{\fr{A}}[\eta]$ is continuous by virtue of, for
example, \cite[Theorem~3.5.6]{tovst1}; the set
$\td{\fr{A}}[\eta](\td{\fr{U}})$ is compact as a continuous image of
a compact. In what follows, is sufficient to use
    $\td{\fr{U}}= cl\, \td{\delta}(\fr{U}).$

    Replacing~$a$ and the compact~$\Xi$
    with the mapping $\{(f,g)\}$ and the compact $\{(0_\mm{X},0\mm{R})\}$,
    we obtain the continuous dependence on~$\eta$
    for the maps $\td{\ph},\td{J}$.
    Now, by virtue of $cl {J}_t(\fr{U})=
    cl \td{J}_t(\td{\delta}\circ\fr{U})=\td{J}_t(\td{\fr{U}})$,
    the condition~${\bf{(e)}}$ holds for $\eta\in\td{\fr{U}}$ too,
     i.e., it holds that
    \begin{equation}
   \label{om}
      \td{J}_t(\eta)\leq J_T(\eta)+\omega(T) \qquad \forall \eta\in\td{\fr{U}},T\in\mm{T},t\in[T,\infty\ra,
\end{equation}
      then,
      $$\limsup_{t\to\infty}\td{J}_t(\eta)\leq J_T(\eta)+\omega(T)\qquad \forall \eta\in\td{\fr{U}},
      T\in\mm{T},$$
passing to the lower limit as $T\to\infty$, we obtain, for arbitrary $\eta\in\td{\fr{U}}$, the existence of the limit $\lim_{t\to\infty}\td{J}_t(\eta)$ (possibly infinite).

Then, for every $t\in\mm{T}$, there exists an $\eta_t\in\td{\fr{U}}$ such that
\begin{equation}
   \label{supu1}
 \fr{R}_t\rav\max_{\eta\in\td{\fr{U}}}\td{J}_t(\eta)=\td{J}_t(\eta_t)
 \qquad
 \forall t\in\mm{T}.
\end{equation}
Since $(\eta_t)_{t\in\mm{T}}$ is in the compact, for the certain  unbounded increasing sequence
  $({t_k})_{k\in\mm{N}}\in\mm{T}$ and the certain $\td{u}^0\in\td{\fr{U}}$,
it is $\eta_{t_k}\to\td{u}^0.$ 
Let us also define
\begin{equation}
   \label{deflim}
   \overline{\fr{R}}\rav\limsup_{t\to\infty}\fr{R}_t,\
   \underline{\fr{R}}\rav\liminf_{t\to\infty}\fr{R}_t,\
  \fr{R}^*\rav
  \sup_{\eta\in\td{\fr{U}}}\lim_{t\to\infty}\td{J}_t(\eta),\
  \fr{R}^0\rav\lim_{t\to\infty} \td{J}_t(\td{u}^0).
\end{equation}
Now,
  $$\fr{R}_{t_k}=\td{J}_{t_k}(\eta_{t_k})\leqref{om}\td{J}_{t_i}(\eta_{t_k})+\omega(t_i)
  \qquad \forall i,k\in\mm{N} (i<k),$$
passing to the upper limit as $k\to\infty$ and then as $i\to\infty$, we obtain
  $\overline{\fr{R}}\leq\td{J}_{t_i}(\td{\fr{u^0}})+\omega(t_i)
  $
and $\overline{\fr{R}}\leq\fr{R}^0$.
Thus, for all $T\in\mm{T}$
\begin{eqnarray*}
 \overline{\fr{R}}\leq\fr{R}^0\ravref{deflim}\lim_{t\to\infty}
\td{J}_t(\td{u}^0)\leq
 \sup_{\eta\in\td{\fr{U}}}\lim_{t\to\infty}\td{J}_t(\eta)\ravref{deflim}
 \fr{R}^*\\ \leqref{om}
\sup_{\eta\in\td{\fr{U}}}\td{J}_T(\eta)+\omega(T)\ravref{supu1}
\fr{R}_T+\omega(T).
\end{eqnarray*}
Passing to the lower limit as $T\to\infty$, we obtain
  $\overline{\fr{R}}\leq\fr{R}^0\leq\fr{R}^*\leq\underline{\fr{R}},$
it remains to note that
  by virtue of $cl {J}_t(\fr{U})=
    cl \td{J}_t(\td{\delta}\circ\fr{U})=\td{J}_t(\td{\fr{U}})$,
it holds that
  $$\lim_{t\to\infty}
  \sup_{u\in\fr{U}}{J}_t(u)=\lim_{t\to\infty}\max_{\eta\in\td{\fr{U}}}\td{J}_t(\eta)=
    \lim_{t\to\infty}\fr{R}_{t}=\fr{R}^0=\fr{R}^*.\qquad\qquad\bo
    $$

\begin{remark}
 As it was shown in~1),
 for each generalized control there exists the  sequence of controls from
$\fr{U}$ that converges (in the topology~$\td{\fr{U}}$) to it.
\end{remark}
\begin{remark}[turnpike property]
\label{serf}
Item~4) actually shows more. It shows that the uniformly  overtaking  optimal control $ \td{u}^0$ can be obtained as a limitary point of the sets $\arg \max_{\eta\in\td{\fr{U}}}
  \td{J}_t(\eta)\in(comp)(\td{\fr{U}})$ as $t\to\infty.$
\end{remark}
\begin{remark}
As it was shown in~4), the limit $$\int_{\mm{T}}\int_{U(t)}g(t,\td{\ph}[\eta](t),u()\eta(du)dt\rav\lim_{T\to\infty}\td{J}_T(\eta) $$
is defined (though it may be infinite) for all $\eta\in\td{\fr{U}}$.
\end{remark}

\medskip

 Note that not only did the paper~\cite{dm} prove the theorem of existence of an optimal solution based on the condition~$\bf{(e)}$ but it also discussed the proof of such theorems based on the inverse limit. To the best of author's knowledge, there is only one paper~\cite{my} besides the previous one in the control theory that explicitly employs the notion of inverse limit.

 There are many existence theorems, for example, \cite{bald}, \cite{car1},\cite{car},\cite{carm}. The results obtained in Proposition~\ref{spectr} have much in common with paper~\cite{carm} (in terms of~\cite{carm}, the obtained~$\td{u}^0$ is strongly optimal).
 Note that if the initial set~$\fr{U}$ does not contain a uniformly  overtaking  optimal    control, we may pass to Gamkrelidze
 controls by increasing the dimension of the set $\mm{U}$ in $m+1$ times.  (For details of such bicompact extension, see~\cite{ga}, \cite{carm}).  These controls also form a compact and the items
1)-3) of Proposition \ref{spectr} hold from them; therefore,  there always exists a uniformly  overtaking optimal control among such finite-dimensional controls.

As a corollary, we assume the uniformly  overtaking optimal control~$u^0$ to exist among the elements of~$\fr{U}$, and
denote the trajectory that corresponds to~$u^0$ by~$x^0$.
We also keep the denotation $\td{u}^0\rav\td{\delta}\circ u^0.$

We are also interested in the degree of closeness of various generalized controls for large~$t$.
 Let $w:\mm{T}\times \mm{U}\to\mm{T}$
  be an integrally bounded Carath\'{e}odory map. For
all $\tau\in\mm{T}$ and $\eta\in\td{\fr{U}}$, let us introduce
  $$\displaystyle \fr{L}_w[\eta](\tau)\rav\int_0^\tau\int_{U(t)} w(t,u)\eta(t)(du)dt.$$
  Let us denote by $(Fin)(u^0)$ the family of $\eta\in\td{\fr{U}}$ such that
 $\eta|_{[T,\infty \rangle}=\td{u}^0|_{[T,\infty\rangle}$
 for the certain $T\in\mm{T}$.
Let us assume that $\fr{L}_w[\td{u}^0]\equiv 0$, and for every
$\eta\in (Fin)(u^0)$ from $\fr{L}_w[\eta](\tau)=0$ for all
$\tau\in\mm{T}$ it follows that $\eta$ equals $\td{u}^0$ a.e. on
$[0,\tau]$. The set of such~$w$ is denoted by $(Null)(u^0)$.

\section{The necessary conditions of optimality}  \label{sec:2}
\ee
\subsection{ Relations of the Maximum Principle} \label{sec:21}
Let  the Hamilton--Pontryagin function
  $\ct{H}:\mm{X}\times Gr\,{U}\times\mm{T}\times\mm{X}\to\mm{R}$
 be given by
 $$  \ct{H}(x,t,u,\lambda,\psi)\rav\psi f\big(t,x,u\big)+\lambda
 g\big(t,x,u\big).$$
Let us introduce the relations
 \begin{subequations}
\begin{eqnarray}
       \dot{x}(t)&=& f\big(t,x(t),u(t)\big);\label{sys_x}
\\
   \label{sys_psi}
       \dot{\psi}(t)\in&-&{\partial_x
       \ct{H}\big(x(t),t,u(t),\lambda,\psi(t)\big)};\\
   \label{maxH}
\sup_{p\in U(t)}\ct{H}\big(x(t),t,p,\lambda,\psi(t)\big)&=&
        \ct{H}\big(x(t),t,u(t),\lambda,\psi(t)\big);\\
   x(0)=0,\ &\ &
       \ \ ||\psi(0)||_\mm{X}+\lambda=1.
   \label{dob}
   \end{eqnarray}
          \end{subequations}

It is easily seen that for each $u\in\fr{U}$, for each initial
condition, system~\rref{sys_x}--\rref{sys_psi} has a local solution,
and each solution of these relations can be extended to the
whole~$\mm{T}.$

Let us denote by~$\fr{Y}$ the family of all solutions
   $(x,u,\lambda,\psi)\in C_{loc}(\mm{T},\mm{X})\times\fr{U}\times[0,1]\times C_{loc}(\mm{T},\mm{X})$
of system~\rref{sys_x}--\rref{sys_psi},\rref{dob} on~$\mm{T}$, and
let us denote by~$\fr{Z}$ the set of solutions from~$\fr{Y}$ for
which~\rref{maxH} also holds a.e. on~$\mm{T}$.

\medskip

Let us introduce such conditions for generalized controls; namely,
under initial condition~\rref{dob} let us consider
  \begin{subequations}
\begin{eqnarray}
   \label{sys_x_}
       \dot{x}(t)&=& \int_{U(t)} f(t,x(t),u) \eta(t)(du);\\
   \label{sys_psi_}
       \dot{\psi}(t)\in&-&\int_{U(t)} {\partial_x
       \ct{H}\big(x(t),t,u,\lambda,\psi(t)\big)}\eta(t)(du);\\
   \label{maxH_}
        \sup_{p\in
        U(t)}\ct{H}\big(x(t),t,p,\lambda,\psi(t)\big)&=&
        \int_{U(t)}\ct{H}\big(x(t),t,u,\lambda,\psi(t)\big)\eta(t)(du).
   \end{eqnarray}
  \end{subequations}
Similarly, for each $\eta\in\td{\fr{U}}$ for each initial condition,
system~\rref{sys_x_}--\rref{sys_psi_} has a local solution that can
be extended to the whole~$\mm{T}$.

Let us denote by~$\td{\fr{Y}}$ the family of all solutions
   $(x,\eta,\lambda,\psi)\in C_{loc}(\mm{T},\mm{X})\times\td{\fr{U}}\times[0,1]\times C_{loc}(\mm{T},\mm{X})$
 of system
   \rref{dob}--\rref{sys_psi_}.
Let us also introduce~$\td{\fr{Z}}$, the family of
$(x,\eta,\lambda,\psi)\in\td{\fr{Y}}$ such that~\rref{maxH_} also
holds a.e. on~$\mm{T}$.

Let us note that for every $\eta\in\td{\fr{U}}$, the family of all
solutions $(x,\eta,\lambda,\psi)\in\td{\fr{Y}}$ of
system~\rref{dob}--\rref{sys_x_} 
 on~$\mm{T}$ for given
control~$\eta$ is compact by virtue of~\cite[Theorem 3.4.2]{tovst1}.
Moreover, this compact-valued map  is upper semicontinuous in~$\eta$.
Indeed, the right-hand side of~\rref{sys_x_}--\rref{sys_psi_} is
convex and integrally bounded, upper semicontinuous in~$\eta$, and it
is measurable for each fixed $x,\psi$; therefore, it has a measurable
selector~(\cite[Lemm
   2.3.11]{tovst1}); moreover, all local solutions of~\rref{sys_x_}--\rref{sys_psi_} can be extended to the whole~$\mm{T}$.
   Since all the conditions of~\cite[Theorem 3.5.6]{tovst1}
   are satisfied, the mapping is upper semicontinuous.
   Therefore,~$\td{\fr{Y}}$ and~$\td{\fr{Z}}$ are compact,
   as the graphs of this mapping on the compact subdomain of its domain.

Note that by \cite[Theorem 2.7.5]{clarke} always holds the inclusion:
\begin{equation}
   \label{psin}
\partial_x \int_{U(t)} {
       \ct{H}\big(x,t,u,\lambda,\psi\big)}\eta(t)(du)\subset
\int_{U(t)} {\partial_x
       \ct{H}\big(x,t,u,\lambda,\psi\big)}\eta(t)(du).
\end{equation}

\subsection{ The necessity of the Maximum Principle} \label{sec:22}

\begin{theorem}
\label{1} Assume conditions  $\bf{(u),(fg)}$. For each
uniformly overtaking  optimal   pair $(x^0,u^0)\in
C(\mm{T},\mm{X})\times\fr{U}$ for problem~\rref{sys}--\rref{opt},
there exist  $\lambda^0\in [0,1],\psi^0\in C(\mm{T},\mm{X})$ such
that the relations of the Maximum Principle \rref{sys_x}--\rref{dob}
hold;
    i.e., $(x^0,u^0,\lambda^0,\psi^0)\in\fr{Z}$.
\end{theorem}
\doc

Let us fix a certain unbounded monotonically increasing sequence
$(\tau_n)_{n\in\mm{N}}\in \mm{T}^\mm{N}$. Let us also consider an
arbitrary sequence $(\gamma_n)_{n\in\mm{N}}\in \mm{T}^\mm{N}$ that
converges to zero with the property $\omega^0(\tau_n)/\gamma_n\to 0$,
 where the function~$\omega^0$ was taken from \rref{555}.
 For example, $\gamma_n\rav \sqrt{\omega^0(\tau_n)}$ will suffice.

 Fix a $w\in(Null)(u^0)$.
 For each $n\in\mm{N}$ let us consider
the problem
$$   J_{\tau_n}(\eta)-\gamma_n
\fr{L}_w[\eta](\tau_n)=\int_0^{\tau_n}\int_{U(t)}
  g(t,\td{\ph}[\eta](t),u)\eta(t)(du)dt-\gamma_n
  \fr{L}_w[\eta](\tau_n)
  \to\max.$$
Here, the functional is bounded from above by the number
$J_{\tau_n}(u^0)+\omega^0(\tau_n)$, therefore, it has the supremum. Every
summand continuously depends on~$\eta$, which covers the
compact~$\td{\fr{U}}$; therefore, there is an optimal solution for
this problem in~$\td{\fr{U}}$; let us denote one of them by
$(x^n,\eta^n)$.

  Let    the function
  $\ct{H}_{\tau_n}:\mm{X}\times Gr\,{U}\times\mm{T}\times\mm{X}\to\mm{R}$
 be given by
 $$\displaystyle \ct{H}_{\tau_n}(x,t,u,\lambda,\psi)\rav\ct{H}(x,t,u,\lambda,\psi)-
  \gamma_n w(t,u).$$
Then, by the Clarke form~\cite[Theorem~5.2.1]{clarke} of the
Pontryagin Maximum Principle, there exists
$(\lambda^n,\psi^n)\in\mm{T}\times C([0,n],\mm{X})$ such that
relation~\rref{dob} and
 the transversality condition at the free endpoint $\psi^n(\tau_n)=0$
 hold,
 and
 \begin{eqnarray}
   \label{sys_max_}
   \sup_{p\in
        U(t)}\ct{H}_{\tau_n}\big(x^n(t),t,p,\lambda^n,\psi^n(t)\big)&=&
           \int_{U(t)}\ct{H}_{\tau_n}\big(x^n(t),t,u,\lambda^n,\psi^n(t)\big)
   \eta^n(t)(du),\\
\dot{\psi}^n(t)&\in& -\partial_x
\int_{U(t)}\ct{H}_{\tau_n}\big(x^n(t),t,u,\lambda^n,\psi^n(t)\big)
   \eta^n(t)(du) \nonumber
 \end{eqnarray}
    also hold for a.a. $t\in[0,\tau_n]$.
 By \rref{psin}, $(x^n,\eta^n,\lambda^n,\psi^n)\in\mm{T}\times C([0,n],\mm{X})$
 satisfy the relations \rref{dob}--\rref{sys_psi_}, \rref{sys_max_} a.~e.
on the $[0,\tau_n].$

Let us extend the  $(x^n,\eta^n,\lambda^n,\psi^n)$ to
$[\tau_n,\infty\rangle$ by the generalized control
$\td{u}^0|_{[\tau_n,\infty\rangle}$. Then, $\eta^n\in (Fin)(u^0)$.
Let us denote by~$\fr{Z}^n$ the set of $(x,u,\lambda,\psi)$ that
satisfy relations
  \rref{dob}--\rref{sys_psi_}
a.~e. on $\mm{T}$, satisfy relation~\rref{sys_max_} a.~e.  on
$[0,\tau_n\rangle$, and possess the property
$\td{u}^0|_{[\tau_n,\infty\rangle}=
   \eta^n|_{[\tau_n,\infty\rangle}$. Now we have $(x^n,\eta^n,\lambda^n,\psi^n)\in\fr{Z}^n$
for every $n\in\mm{N}$.

Let us note that all~$\fr{Z}^n$ are closed and, since these sets are
contained
 in the compact $\td{\fr{Y}}$, these sets are also compact. Hence, the sequence
  $(x^n,\eta^n,\lambda^n,\psi^n)_{n\in\mm{N}}$ has the limit point
  $(x^{00},\eta^0,\lambda^0,\psi^0)\in\td{\fr{Y}}$.
  Passing, if  necessary,   to a subsequence,  we may  assume that
  it is the limit of the sequence itself.

For a fixed~$x$, the set of $u\in U(t)$ that realize the maximum
in~\rref{sys_max_} has a measurable selector by virtue
of~\cite[Theorem 3.7]{select}. By \cite[Lemm 2.3.11]{tovst1}, it
exists if we put an arbitrary continuous function~$x$ into~$\ct{H}$.
Besides, since relation~\rref{sys_max_} also depends on $x,\psi$ and
on the parameters~$\gamma$ and~$\lambda$ upper semicontinuously, and
all the relations are integrally bounded on bounded sets;
 by virtue of~\cite[Theorem 3.5.6]{tovst1}, on each finite interval
 for the funnels of solutions of~\rref{sys_x}--\rref{sys_psi} that
 satisfy~\rref{sys_max_}, we have upper semicontinuity by
 $\gamma,\lambda$.
  In particular, for $\gamma_n\to 0$, $\lambda^n\to \lambda^0$,
  we obtain the fact that the upper limit of the compacts $\fr{Z}^n$
   is included in $\td{\fr{Z}}$. Hence,
   $(x^{00},\eta^0,\lambda^0,\psi^0)\in\td{\fr{Z}}$.

 On the other side, by $w\in(Null)(u^0)$ and by optimality of~$\eta^n$, $u^0$ for their problems, we obtain
$$ \td{J}_{\tau_n}(\eta^n)-\gamma_n \fr{L}_w[\eta^n](\tau_n)\geq
 J_{\tau_n}(u^0)\geqref{555} \td{J}_{\tau_n}(\eta^n)-\omega^0(\tau_n)
 $$
therefore, we have
   $\displaystyle
      \gamma_n \fr{L}_w[\eta^n](\tau_n)
 \leq \omega^0(\tau_n).$
 By virtue of
    $\td{u}^0|_{[\tau_n,\infty\rangle}=
   \eta^n|_{[\tau_n,\infty\rangle},$
we obtain
\begin{equation}
   \label{to_w}
   \fr{L}_w[\eta^n](\tau)\leq \omega^0(\tau_n)/\gamma_n
   \qquad\forall \tau\in\mm{T}.
\end{equation}
For each $\tau\in\mm{T}$, passing to the limit as $n\to\infty$, we
obtain that $\fr{L}_w[\eta^0]\leq 0$; i.e.,
$\fr{L}_w[\eta^0](\tau)=0$ for all $\tau\in\mm{T}$. Since
$w\in(Null)(u^0)$, we have $\eta^0=\td{u}^0$ a.e. on $\mm{T}$, hence
   $x^{00}=x^0$ and $(x^0,u^0,\lambda^0,\psi^0)\in{\fr{Z}}$.
    Moreover, from~\rref{to_w}, we have $||\fr{L}_w[\eta^n]||_C\to 0$.\bo

\medskip
We have additionally proved that
\begin{remark}
 \label{1014}
Under conditions $\bf{(u),(fg)}$, for each optimal pair      $(x^0,u^0)\in\fr{X}\times\fr{U}$
for problem~\rref{opt}, for each weight $w\in(Null)(u^0),$
 for each unbounded increasing sequence      $(\tau_n)_{n\in\mm{N}}\in \mm{T}^\mm{N}$, we have constructed the sequence
  $(x^n,\eta^n,\lambda^n,\psi^n)_{n\in\mm{N}}\in \td{\fr{Y}}^\mm{N}$ that possesses the following properties:
  \begin{description}
    \item[1)] This sequence (as a sequence from
  $C_{loc}(\mm{T},\mm{X})\times \td{\fr{U}}\times\mm{T}\times C_{loc}(\mm{T},\mm{X})$)
  converges to the certain $(x^0,\td{u}^0,\lambda^0,\psi^0)\in
   {\fr{Z}}$;
    \item[2)] $||\fr{L}_w(\eta^n)||_C\to 0$;
    \item[3)] $\td{J}_{t_n}(\eta)-{J}_{t_n}(u^0)\to 0$, and $\psi^n(t_n)=0$ for each
$n\in\mm{N},$ where $(t_n)_{n\in\mm{N}}$ is a certain subsequence  of
  $(\tau_n)_{n\in\mm{N}}\in \mm{T}^\mm{N}.$
  \end{description}
\end{remark}

\subsection{The simplest condition of transversality}
\label{sec:23}

 However, the relations of the Maximum Principle are
incomplete, since \rref{sys_x}--\rref{dob}
do not contain a condition on the right endpoint. There are several
variants of such additional conditions (for details,
see~\cite[Sect.~6,12]{kr_as},\cite{ss}); in this paper we investigate
the modifications of the condition
\begin{subequations}
\begin{equation}
   \label{trans}
       \displaystyle \lim_{t\to\infty} \psi(t)=0.
   \end{equation}

Let us formulate the propositions in terms of the stability of~$\psi$
such that a condition  would be necessary.

\medskip

  {\it \bf Condition}~{($\psi$)}: There exists a weight $w\in(Null)(u^0)$
 such that for every solution
 $(x^0,u^0,\lambda^0,\psi^0)\in{\fr{Z}}$,
 the Lagrange multiplier~$\psi^0$ is stable under
 $\fr{L}_w-$small perturbations of system
 \rref{sys_x}--\rref{sys_psi};\ 
i.e., for every $\epsi\in\mm{R}_{> 0}$, there exist a number
$\delta\in\mm{R}_{> 0}$ and a neighborhood $\Upsilon\subset
C_{loc}(\mm{T},\mm{X})\times \td{\fr{U}}\times[0,1]\times
C_{loc}(\mm{T},\mm{X})$ of the solution
$(x^0,\td{u}^0,\lambda^0,\psi^0)$
 such that for every solution
$(x,\eta,\lambda,\psi)\in\Upsilon\cap\td{\fr{Y}}$ from
$||\fr{L}_w[\eta]||_C<\delta$ it follows that
$||\psi-\psi^0||_C<\epsi$.

\medskip

\begin{proposition}
\label{3}
 Assume conditions $\bf{(u),(fg)}$ hold. For each uniformly  overtaking  optimal
pair
    $(x^0,u^0)\in C(\mm{T},\mm{X})\times\fr{U}$
satisfying $\bf{(\psi)}$, for each unbounded increasing sequence,
$(\tau_n)_{n\in\mm{N}}\in \mm{T}^\mm{N}$
 there exists
 $(x^0,u^0,\lambda^0,\psi^0)\in\fr{Z}$ such that
    \begin{equation}
   \label{partlim}
       \liminf_{n\to\infty}||\psi^0(\tau_n)||_{\mm{X}}=0
   \end{equation}
holds.
\end{proposition}
          \end{subequations}
\doc
Let us choose the certain $\epsi\in\mm{R}_{> 0}$, and let us take
$\Upsilon\subset\td{\fr{Y}}$ and $\delta\in\mm{R}_{> 0}$ from
condition~{($\psi$)}; by Remark~\ref{1014},
 there exists $N\in\mm{N}$ such that for
$n\in\mm{N}, n>N$, it is
  $(x^n,\eta^n,\lambda^n,\psi^n)\in\Upsilon$,
   $||\fr{L}_w[\eta^n]||_C<\delta$; now,
   condition~{($\psi$)} also yields
   $||\psi^n(\tau_n)-\psi^0(\tau_n)||_\mm{X}<\epsi$;
   but  $\psi^n(\tau_n)=0$;
   whence $||\psi^0(\tau_n)||_\mm{X}<\epsi$ for all $n\in\mm{N}, n>N$.
   Since $\epsi\in\mm{R}_{> 0}$ was arbitrary, we have shown \rref{partlim}.
\bo

Note that by linearity of~\rref{sys_psi}, the stability of the variable~$\psi$ implies its boundedness. Therefore, the proved proposition is useless for unbounded adjoint variable~$\psi$.

 Note that, as it follows from \cite[Example~5.1]{stern}, for a uniformly overtaking optimal control, there can be no
 $(x^0,u^0,\lambda^0,\psi^0)\in\fr{Z}$ that satisfies stronger condition~\rref{trans} instead of~\rref{partlim}. On the other side,
 \begin{remark}~\label{partlimz}
 Assume the functions $L_K^f,L_K^g$ are
independent of a compact~$K$, and the mapping
 $T\mapsto L_K^g(T) e^{\int_{[0,T]}L_K^f(t)dt}$ is
 summable  
 on  $\mm{T}$
    (\cite[Hypotesis 3.1 (iv)]{sagara}); therefore, the total variation of~$\psi$
    is a fortiori bounded. Then, $\bf{(\psi)}$ holds and, moreover,
    \rref{partlim}  implies \rref{trans}.
\end{remark}

The even more strong conditions used for proving the Maximum
Principle can be seen, for example,
 in~\cite[(A3)]{Ye} (the Lipshitz
constants were required to decrease exponentially with time). Naturally, the propositions proved there for the condition are
also covered by proposition~\ref{3}.

 One of the most general conditions on~\rref{trans} was shown
in~\cite{norv}. For a control problem without phase restrictions,
the transversality condition from \cite[Theorem 6.1]{norv} follows
from Proposition~\ref{3} and~\cite[Lemm 3.1]{norv}, or from
Remark~\ref{partlimz} and condition \cite[(C3)]{norv}. The Remark
\ref{1014}  automatically yields~\cite[Theorem 8.1]{norv}.

\section{The necessity and the stability} \label{sec:3}
\ee
 The objective at hand is to choose the weight
$w^0\in(Null)(u^0)$ such that condition~{($\psi$)} would follow from a
variety of (nonasymptotic) Lyapunov stability of~$\psi$.

\subsection{On weight $w^0$} \label{sec:32}

 Assume conditions $\bf{(u),(a)}$ hold.
 In what follows, assume $\Xi\rav E$.
 Then, for every position
$(\tau^*,y^*)\in\mm{T}\times{E}$ there exists the unique
solution~$y^0$ of the equation
\begin{equation}
   \label{1667}
    \dot{y}=a(t,y(t),u^0(t)),\quad y(\tau^*)=y^*, \tau^*\in \mm{T}
\end{equation}
that can be extended to the whole time interval $\mm{T}$. It (as an
element of $\td{\fr{A}}(\td{u}^0)\subset C_{loc}(\mm{T},E)$)
continuously depends on $(\tau^*,y^*)\in\mm{T}\times E$. Let us
denote its initial position~$y^0(0)$ by $\beg(\tau^*,y^*)$.

\medskip

\begin{proposition}
\label{dop} Assume $\bf{(u),(a)}$ hold. Let the compact-valued map
$G:\mm{T}\rightsquigarrow E$ be bounded on each compact set, and let
$Gr\, G$ be closed.

Then, there exists $w^0\in(Null)(u^0),$ 
such that for arbitrary $\eta\in\td{\fr{U}}$, $T\in\mm{T}$ for every
$y\in\td{\fr{A}}[\eta]$ from $Gr\, y|_{[0,T]}\subset Gr\, G$ it
follows that
 $$  ||\beg(\tau,y(\tau))-y(0)||_E\leq \fr{L}_{w^0}[\eta](\tau)\qquad\forall \tau\in
 [0,T].$$
 \end{proposition}
\doc

Fix an $n\in\mm{N}$. By continuability, for each $(\tau^*,y^*)\in
Gr\,G|_{[0,n]}$, there exists the position $\beg(\tau^*,y^*);$ by
virtue of the theorem of continuous dependence on initial conditions,
this mapping is continuous; hence, the image
      $$\bar{G}_n\rav\Big\{ e\in Gr\,\,
      y|_{[0,n]}\,\Big|\,\forall y\in\td{\fr{A}}[\td{u}^0],(\tau^*,y(\tau^*))\in
       Gr\,\, G|_{[0,n] } \Big\}$$ is closed;
       by the continuability, this set is bounded and, therefore, compact.
        Therefore, on this set, the function $a(t,y,u^0(t))$
       is Lipshitz continuous with respect to~$y$
       for the certain Lipshitz constant
       $L_n\rav L^a_{\bar{G}_n}\in \ct{L}^1_{loc}(\mm{T},\mm{T})$.
       For all $t\in[0,n]$, define
       $M_n(t)\rav\int_{[0,t]}L_n(\tau)d\tau$.
       Note that this function is absolutely continuous and monotonically nondecreasing.

Fix $n\in\mm{N}$; for all $t\in [n-1,n\ra, u\in \mm{U}$, let us
consider a number
 $$    R(t,u)\rav \sup_{ y\in \bar{G}_n}
 \big|\big|a(t,y,u)-a(t,y,u^0(t)\big)\big|\big|_E.
     $$
 Note that the norm inside is a mapping that is continuous with respect to~$y$ and~$u$,
 and~$y$ assumes values from the compact set; now, for every $u\in \mm{U}$
 by~\cite[Theorem 3.7]{select} the supremum reaches the maximum
for the certain function $y_{max}[u]\in
\ct{L}^1([n,n-1\ra,\bar{G}_n)$. Hence, $R(t,u)$ is measurable with
respect to~$t$ for each $u\in\mm{U}$.

Fix a $t\in [n-1,n\ra$; for each sufficiently small neighborhood $\Upsilon\subset U(t)$, by continuity of
$a(t,\cdot,\cdot)$ on compact $\bar{G}_n\times cl \Upsilon$, there
exists a function $\omega^t\in\Omega$, for which
\begin{equation}
  \label{1111}
   \Big|\,\big|\big|a(t,y,u')-a(t,y,u^0(t)\big)\big|\big|-\big|\big|a(t,y,u'')-a\big(t,y,u^0(t)\big)\big|\big|\,\Big|<
   \omega^t\big({
   \frac{1}{||u'\!-\!u''||}}\big)
\end{equation}
holds for every $y\in\bar{G}_n, u',u''\in \Upsilon\, (u'\neq u'')$.
Without loss of generality,
 assume $R(t,u')\leq R(t,u'')$. 
Now, by definition,
  $R(t,u')\geq
  \big|\big|a(t,y,u')-a\big(t,y,u^0(t)\big)\big|\big|$,
   and, substituting $y\rav y_{max}(u'')(t)$ into~\rref{1111}, we obtain
  $0\leq R(t,u'')-R(t,u')\leq
  \omega^t(1/||u'-u''||);$
  i.e.,~$R$ is continuous   with respect to the variable~$u$ on
each sufficiently small neighborhood $\Upsilon\subset U(t)$;
therefore on~$U(t)$ and $Gr U|_{[n-1,n\ra}$ too. Thus, the function
$R: Gr U|_{[n-1,n\ra}\to\mm{T}$ is a Carath\'{e}odory function.

Let us note that by considering all $n\in\mm{N}$, we define the
Carath\'{e}odory function~$R$ on the whole $Gr U$. Moreover, by
construction, $R(t,u^0(t))\equiv 0$. Hence, it is correct to define
$w^0\in (Null)(u^0)$ by the rule
 $$ w^0(t,u)\rav||u-u^0(t)||+e^{M_n(t)}R(t,u) \qquad \forall n\in\mm{N},(t,u)\in Gr\, U|_{[n-1,n\ra}.$$

Consider arbitrary $n\in\mm{N}$, $\tau^*\in [0,n]$, and
$(\tau^*,y^*_1),(\tau,y^*_2)\in
        \bar{G}_n$. 
For the solutions $y_1,y_2\in \td{\fr{A}}[\td{u}^0]$ of
equation~\rref{1667}, for the initial conditions $y_i(\tau^*)=y^*_i$,
we have        $Gr\, y_i|_{[0,n]}\subset \bar{G}_n$.  Let us
introduce functions
        $$r(t)\rav y_1(t)-y_2(t),\quad
                 W_+(t)\rav e^{M_n(t)}||r(t)||_E\qquad \forall t\in [0,n].$$
By Lipshitz continuity of the right-hand side of~\rref{1667} we
obtain $ ||\dot{r}(t)||_E\geq -L_n(t)
       ||r(t)||_E,$ and
       $$\frac{dW^2_+(t)}{dt}=2L_n(t)W^2_+(t)
       +2e^{2M_n(t)}r(t)\dot{r}(t)
       \geq2L_n(t)W^2_+(t)
       -
       2L_n(t)W^2_+(t)
       =0.
       $$
       Thus, the function~$W_+$ is nondecreasing,
       and finally for all $(\tau,y^*_1),(\tau,y^*_2)\in
        \bar{G}_n$ we have
\begin{equation}
  \label{1037_}
     ||\beg(\tau,y^*_1)-\beg(\tau,y^*_2)||_E={W_+(0)}\leq {W_+(\tau)}=
      e^{M_n(\tau)}||y^*_1-y^*_2||_E.
\end{equation}

 Assume the
 $\eta\in\td{\fr{U}},$ $y\in\td{\fr{A}}[\eta],$
         $T\in\mm{T}$ satisfy $Gr\, y|_{[0,T]}\subset Gr\, G$.
         Fix arbitrary $n\in\mm{N}$ and $\tau_1,\tau_2\in[0,T]\cap[n-1,n\ra,$
        $\tau_1<\tau_2$. There exists the solution $y^0\in\td{\fr{A}}[\td{u}^0]$ that
         satisfies the condition
        $y^0(\tau_1)=y(\tau_1)$; let us also define
        $$r\rav y^0(t)-y(t),\quad W_-(t)\rav e^{-M_n(t)}||r(t)||_E\qquad \forall t\in [\tau_1,\tau_2].$$
By construction of $\bar{G}_n$, we have
        $Gr\, y|_{[\tau_1,\tau_2]},Gr y^0|_{[\tau_1,\tau_2]}\subset \bar{G}_n.$
Now,
 \begin{eqnarray*}
  \frac{dW_-^2(t)}{dt}=
       2e^{-2M_n(t)}r(t)\dot{r}(t)-2L_n(t)W_-^2(t)
       =\\
       2e^{-2M_n(t)}{r}(t)\big(\dot{y}^0(t)\!-\!a(t,{y}(t),u^0(t))+
       a(t,{y}(t),u^0(t))\!-\!\dot{y}(t)\big)-2L_n(t)W_-^2(t)
       \leq \\
          2e^{-2M_n(t)}||r(t)||_E
          \int_{U(t)}R(t,u)
          \eta(t)(du)+2L_n(t)W_-^2(t)-2L_n(t)W_-^2(t)\leq \\
   2e^{-M_n(t)}W_-(t)
          \int_{U(t)}R(t,u)
          \eta(t)(du)\leq 2e^{-2M_n(t)}W_-(t)\frac{d
          \fr{L}_{w^0}[\eta](t)}{dt}.
%
 \end{eqnarray*}
     Since function~$W_-$ is nonnegative,
     for a.~a. $t\in \{t\in[\tau_1,\tau_2]\,|\,W_-(t)\neq 0\}$ we obtain
 \begin{equation}
  \label{1623}
     \frac{dW_-(t)}{dt}\leq
     e^{-2M_n(t)}\frac{d \fr{L}_{w^0}[\eta](t)}{dt}\leq e^{-2M_n(\tau_1)}\frac{d
     \fr{L}_{w^0}[\eta](t)}{dt}.
\end{equation}
This inequality is trivial for
     $[\tau_1,\tau_2]\ni t<\sup \{t\in[\tau_1,\tau_2]\,|\,W_-(t)=0\};$
whence,
\begin{eqnarray*}
||\beg(\tau_2,y^0(\tau_2))-\beg(\tau_2,y(\tau_2))||_E&\leqref{1037_}&
e^{M_n(\tau_2)}||y^0(\tau_2)-y(\tau_2)||_E=     \\
     e^{2M_n(\tau_2)}W_-(\tau_2)&\leqref{1623}&
e^{2M_n(\tau_2)-2M_n(\tau_1)}
          \big(\fr{L}_{w^0}[\eta](\tau_2)-\fr{L}_{w^0}[\eta](\tau_1) \big).
\end{eqnarray*}
But
$\beg(\tau_2,y^0(\tau_2))=y^0(0)=\beg(\tau_1,y^0(\tau_1))=\beg(\tau_1,y(\tau_1)),$
hence, we have
 \begin{equation}
  \label{1649}
     ||\beg(\tau_2,y(\tau_2))-\beg(\tau_1,y(\tau_1))||_E\leq
    e^{2M_n(\tau_2)-2M_n(\tau_1)}
          \big(\fr{L}_{w^0}[\eta](\tau_2)-\fr{L}_{w^0}[\eta](\tau_1) \big).
\end{equation}

Fix arbitrary $t\in [0,T]$. For each $\epsi\in\mm{R}_{>0}$ we can
split interval $[0,t\ra$ into the intervals of the form
$[\tau',\tau''\ra$ such that
 $M_n(\tau'')-M_n(\tau')=\int_{[\tau',\tau''\ra} L_n(t)dt<\epsi$
 and
 $[\tau',\tau''\ra\subset[n-1,n\ra$ for the certain $n\in\mm{N}.$
But,~\rref{1649} holds for every interval, i.e.,
$$||\beg(\tau'',y(\tau''))-\beg(\tau',y(\tau'))||_E\leq
    e^{2\epsi}
          \big(\fr{L}_{w^0}[\eta](\tau'')-\fr{L}_{w^0}[\eta](\tau') \big).$$
Summing for all intervals, by $\beg(0,y(0))=y(0)$ and by the triangle
inequality, we obtain
    $||\beg(t,y(t))-y(0)||_E\leq
    e^{2\epsi}
          \fr{L}_{w^0}[\eta](t)$
          for every $t\in [0,T]$. Arbitrariness of $\epsi\in\mm{R}_{>0}$
           completes the proof of the proposition.
 \bo

\subsection{The partial Lyapunov stability} \label{sec:33}

Assume~$E$ can be represented in the form $E=E_p\times E_q$ for some
finite-dimensional Euclidean subspaces $E_p$ and $E_q$. Let us denote
the projections of the map~$a$ to the subspaces~$E_p$ and~$E_q$
by~$b$ and~$c$, respectively. Now, the system~\rref{a} can be written
in the form
\begin{equation}
   \label{aa}
   \dot{p}= b(t,p,q,u),\quad \dot{q}= c(t,p,q,u),\qquad
 (p,q)(0)=\xi\in E,\quad u\in U(t);
\end{equation}
Then, it is possible to say that
 for all $\eta\in\td{\fr{U}}$, the set $\td{\fr{A}}[\eta]$ contains
pairs of functions $(p,q)\in C_{loc}(\mm{T},E_p)\times
C_{loc}(\mm{T},E_q)$. For every $\xi\in E$, let us denote by
$y^0_\xi\rav(p^0_\xi,q^0_\xi)\in\td{\fr{A}}[\td{u}^0]$ the unique
solution of~\rref{1667} for $\tau^*=0,y^*=\xi.$

\begin{definition}
Consider a closed set $G_0\subset E$ and $\xi\in G_0$. We say that {\it
the solution~$y^0_\xi$ of equation~\rref{1667} has Lyapunov stable
component~$p^0_\xi$  in domain~$G_0$} if for each
$\epsi\in\mm{R}_{>0}$ there exists $\delta(\epsi,y)\in\mm{R}_{>0}$
such that for each
    $\xi'\in G_0$
 from  $||\xi'-\xi||_{E}<\delta(\epsi,y)$ it follows
  that $||p^{0}_{\xi'}(s)-p^{0}_\xi(s)||_E<\epsi$ for all $s\in\mm{T}$.
\end{definition}

\begin{proposition}
\label{zvez} Assume~$\bf{(u),(a)}$ holds. Suppose there is a closed
set $G_0\subset E$ and a compact $K_0\in(comp)(G_0)$ such that for
each $\xi\in K_0$ the solution $y^0_\xi$ of equation~\rref{1667}
has Lyapunov stable component~$p^0_\xi$ in~$G_0$.

Then, for each $\epsi\in\mm{R}_{>0}$, there exists a number
 $\delta\in\mm{R}_{>0}$ 
such that for all $\eta\in \td{\fr{U}},
        y=(p,q)\in\td{\fr{A}}[\eta]$
from $y(0)\in K_0, ||L_{w^0}[\eta]||_C<\delta$, and $\beg(t,y(t))\in
G_0$ for all $t\in \mm{T},$ it follows that
        $||p-p^0_{y(0)}||_C<\epsi.$
\end{proposition}

 \doc
Consider a compact $K_>\rav\{\xi\in G_0\,|\,\exists \xi_0\in K_0\,\,
||\xi-\xi_0||_E\leq 1\}$. To each $t\in\mm{T}$, let us assign the set
$G(t)\rav\{y(t)\,|\,\eta\in\td{\fr{U}},y\in\td{\fr{A}}[\eta],y(0)\in
K_>\}$. The obtained map $G$ is compact-valued and continuous; in
particular, its graph is closed. Now we can use Proposition~\ref{dop}
for the multi-valued map~$G$ and fix the weight $w^0\in(Null)(u^0)$
which exists by this Proposition.

Define
  $$M(\xi',\xi'')\rav\sup_{t\in\mm{T}}
||p^{0}_{\xi'}(t)-p^{0}_{\xi''}(t)||_{E_p}\in\mm{T}\cup\{+\infty\}\qquad
\forall \xi',\xi''\in K_{>}.$$
 For all $\xi\in K_{0}$,
the stability of the component~$p^0_\xi$ implies that the map~$M$ is
finite and continuous at the point $(\xi,\xi)\in K_{>} \times
K_{>}.$

Fix an $\epsi\in\mm{R}_{>0}$; choose for every $\xi\in  K_{0}$ its
$\delta(\epsi/2,y^{0}_\xi)\in\la 0,1/2]$; now, we have also chosen
the $\delta(\epsi/2,y^{0}_\xi)-$neighborhood of the point $(\xi,\xi)$
  (in $ K_>\times K_>$).
From the obtained cover of the diagonal~$\Delta$ of the set
$K_0\times K_0$, let us select a finite subcover; it induces certain
open neighborhood~$\Upsilon$ of the diagonal~$\Delta$.
Let~$\delta(K_0)$ be the minimum distance from the diagonal~$\Delta$
to the boundary of the neighborhood~$\Upsilon$. Now, for all $\xi'\in
K_>$, $\xi\in{K}_0$ from $||\xi'-\xi||_E<\delta(K_0)$ it follows that
$(\xi',\xi)\in\Upsilon$; i.e., for some $\xi''\in K_0$ we have
$M(\xi,\xi''),M(\xi'',\xi')<\epsi/2$, whence $M(\xi,\xi')<\epsi$.
Thus,
\begin{equation}
\label{1291} (||\xi'-\xi||_E<\delta(K_0))\quad\Rightarrow\quad
  (||p^{0}_{\xi'}-p^{0}_\xi)||_C<\epsi)\qquad\forall \xi\in{K}_0,\xi'\in K_>.
\end{equation}

 Suppose the $\eta\in\td{\fr{U}},
y=(p,q)\in\td{\fr{A}}[\eta]$ satisfy
$L_{w^0}[\eta](t)<\delta,\xi_1(t)\rav \beg(t,y(t))\in G_0$ for all
$t\in\mm{T}$. For $K_0\subset K_>=G(0)$, the definition
  $$T_0\rav\sup \{T\in\mm{T}\,|\, Gr\, \xi_1|_{[0,t]}\subset K_>\qquad \forall t\in[0,T\ra\}
  \in \mm{T}\cup\{+\infty\} $$
is correct, although~$T_0$ can be infinite. Hence,  we have $Gr\,
{y}|_{[0,t]} \subset Gr\, G$ for all $t\in [0,T\ra$. Now, from
Proposition~\ref{dop}, we obtain
\begin{equation}
\label{1364}
 ||{\xi}_1(t)-y(0)||_E=||\beg(t,{y}(t))-{y}(0)||_E\leq
 \fr{L}_{w^0}[\eta](t)<\delta(K_0)\qquad \forall t\in [0,T_0\ra.
\end{equation}
For every $t\in [0,T_0\ra$, let us substitute
$\xi=y(0),\xi'\rav{\xi}_1(t)\in K_>$ in~\rref{1291}; from the
equality ${p}^{0}_{\xi_1(t)}(t)={p}(t)$ we obtain
$||{p}(t)-{p}^{0}_{y(0)}(t)||_E<\epsi$ for all $[0,T_0\ra$. To
conclude the proof, it remains to prove that $T_0=\infty.$

Suppose $T_0\in\mm{T}$; by construction of~$T_0$, for each $\tau\in
\la T_0,\infty\ra$, we have
   $Gr\, \xi_1|_{\la T_0,\tau]}\not\subset Gr\, K_>;$ but $Gr\, \xi_1\subset G_0$.
    Then,
   $\rho(\xi_1(T_0),G_0\setminus K_>)=0$, and, in particular,
   by construction of~$K_>$, we have $||\xi_1(T_0)-y(0)||_E\geq 1$.
   However, passing to the limit in~\rref{1364} yields
    $||\xi_1(T_0)-y(0)||_E\leq\delta(K_0)\leq 1/2$. The acquired contradiction proves that
   $T_0=\infty$. \bo

\subsection{The necessity of the transversality
condition~\rref{partlim}} \label{sec:34}

Everywhere further, we assume the following condition holds:

 {\it \bf Condition~}${\bf{(\partial)}}:$ for the maps
      $(t,x)\in\mm{T}\times\mm{X}\times\mm{U}\to f(t,x,u)\in\mm{X}$ and
      $(t,x)\in\mm{T}\times\mm{X}\times\mm{U}\to g(t,x,u)\in\mm{R}$
on their respective domains, there exist partial derivatives in~$x$
that are integrally bounded (on each compact) locally Lipshitz
continuous Carath\'{e}odory maps.
\medskip

Under this condition, the set ${\partial_x
\ct{H}(x(t),t,u^0(t),\lambda,\psi(t))}$ is also a single-element set,
 therefore system~\rref{sys_x}--\rref{sys_psi} can be
rewritten for $u=u^0$ in the form
  \begin{subequations}
  \begin{eqnarray}
   \label{sys_psii}
       \dot{\psi}(t)=-\frac{\partial\ct{H}}{\partial
       x}(x(t),t,u^0(t),\lambda,\psi(t)),\\
   \label{sys_xxl}
       \dot{x}(t)=f(t,x(t),u^0(t)),\\
   \label{sys_lam}
       \dot{\lambda}=0.
   \end{eqnarray}
   \end{subequations}

\begin{corollary}
\label{s3}
 Assume conditions $\bf{(u),(fg),(\partial)}$ hold.
 Let the pair $(x^0,u^0)\in C_{loc}(\mm{T},\mm{X})\times\fr{U}$
  be  uniformly  overtaking  optimal   for problem~\rref{sys}--\rref{opt}.
   If for each  solution $(\psi^0,x^0,\lambda^0)$
 of system
    \rref{sys_psii}--\rref{sys_lam}
    with initial conditions from
$K_0\rav\mm{D}\times[0,1]\times\{0_\mm{X}\}$
the component $\psi^0$ is partially Lyapunov stable in
$G_0\rav\mm{X}\times[0,1]\times\mm{X}$.

 Then, the result of
Proposition~\ref{3} holds.
\end{corollary}

  \doc
In~\rref{aa}, it is sufficient to define $E_p\rav\mm{X},$
$E_q\rav\mm{X}\times\mm{R}$ and to take $\psi$ and $(x,\lambda)$ for
$p$ and $q=(q_1,q_2)$, and to take for~$b$ and~$c$ the right-hand
sides of~\rref{sys_psii} and~\rref{sys_xxl}--\rref{sys_lam},
respectively. Now\
Proposition~\ref{zvez} guarantees {($\psi$)}, i.e., all conditions of
Proposition~\ref{3} are met. This proves Corollary \ref{s3}. \bo

For~$G_0$, we can take the image
    $\{\beg(t,\psi(t),x(t),\lambda)\,|\,
         (x,u,\lambda,\psi)\in\td{\fr{Y}},t\in\mm{T}\}.$

Using Proposition~\ref{dop} and Remark~\ref{1014} instead of Proposition~\ref{zvez} and Corollary~\ref{s3} in this proof, we obtain
\begin{remark}
 \label{1514}
Under conditions $\bf{(u),(fg),(\partial)}$, for each uniformly
overtaking  optimal  pair $(x^0,u^0)\in
C_{loc}(\mm{T},\mm{X})\times\fr{U}$ for problem~\rref{opt}, for each
unbounded increasing sequence $(\tau_n)_{n\in\mm{N}}\in
\mm{T}^\mm{N}$, we have constructed the sequence
  $(x^n,\eta^n,\lambda^n,\psi^n)_{n\in\mm{N}}\in \td{\fr{Y}}^\mm{N}$ such that
  \begin{description}
    \item[1)] this sequence (as a sequence from
  $C_{loc}(\mm{T},\mm{X})\times \td{\fr{U}}\times\mm{T}\times C_{loc}(\mm{T},\mm{X})$) converges to a certain $(x^0,\td{u}^0,\lambda^0,\psi^0)\in
   {\fr{Z}}$;
\item[2)] the graphs
   $Gr (x^n,\lambda^n,\psi^n)$ of its elements are contained within the thinning funnels of solutions of
    system~\rref{sys_psii}--\rref{sys_lam}; i.e., for a sequence $(\delta_n)_{n\in\mm{N}}\in\mm{R}_{\geq 0}^\mm{N}$ that tends to $0$, we have
   $$\varkappa(t,(\psi^n,x^n,\lambda^n))\in (\psi^0(0),0,\lambda^0)+
   \delta_n\mm{D}\times\delta_n\mm{D}\times[-\delta_n,\delta_n]
    \qquad \forall t\in\mm{T},n\in\mm{N}.$$
    \item[3)] $\td{J}_{t_n}(\eta)-{J}_{t_n}(u^0)\to 0$, and $\psi^n(t_n)=0$ for each
$n\in\mm{N},$ where $(t_n)_{n\in\mm{N}}$ is a certain subsequence  of
  $(\tau_n)_{n\in\mm{N}}\in \mm{T}^\mm{N}.$
  \end{description}

\end{remark}

\subsection{Modifications of transversality condition~\rref{partlim}}
\label{sec:35}

 In certain cases, if the Lagrange multiplier~$\psi$ is not
stable, but  we know that certain components of the vector
variable~$\psi$ are stable, or we know the rate of its growth. Then
we can try to select the mapping $A_*:\mm{T}\to\mm{L}$, may help to
modify condition~\rref{partlim}, and use the condition
\begin{subequations}
     \begin{equation}
   \label{partlim_1}
       \liminf_{t\to\infty}
       ||\psi^0(t)A_*(t)||_\mm{X}=0
   \end{equation}
 for certain map~$A_*:\mm{T}\mapsto\mm{L}.$

Here are the examples of such maps $A_*$:
 one that maps the unity matrix $A_*(t)\equiv1_\mm{L};$
 some ``scalar'' multiplier $A_*(t)\equiv r(t)1_{\mm{L}};$
 a mapping $A_*(\cdot)\equiv D$ with the diagonal matrix~$D$;
 the condition $\psi(t)x(t)\to 0$, which is often used as the
 sufficient condition, can also be reduced to this form.

Let us assume that for all $\eta\in\td{\fr{U}},\xi\in\mm{X}$, we
have chosen the measurable mapping
  $A^\eta_\xi:\mm{T}\to\mm{L}$. Assume $A^\eta_\xi(0)=1_\mm{L}$
  for all $\eta\in\td{\fr{U}},\xi\in\mm{X}.$
  Define $A_*\rav A^{\td{u}^0}_{0_\mm{X}}.$

\medskip

  {\it \bf Condition}~{($\psi A$)}: 
    There exists a weight $w\in(Null)(u^0)$ such that for every solution
 $(x^0,u^0,\lambda^0,\psi^0)\in{\fr{Z}}$
for each $\epsi\in\mm{R}_{> 0}$ there exist a number
$\delta\in\mm{R}_{> 0}$ and a neighborhood $\Upsilon\subset
C_{loc}(\mm{T},\mm{X})\times \td{\fr{U}}\times[0,1]\times
C_{loc}(\mm{T},\mm{X})$ of the solution
$(x^0,\td{u}^0,\lambda^0,\psi^0)$
   such that for every solutions $z\rav(x,\eta,\lambda,\psi)\in\Upsilon\cap \td{\fr{Y}}$,
    from $||\fr{L}_w[\eta]||_C<\delta$ it follows that $||\psi A^\eta_{x(0)}-
    \psi^0 A^{\td{u}^0}_{0_\mm{X}}||_C<\epsi$.

\medskip

\begin{proposition}
\label{3_} Assume conditions $\bf{(u),(fg)}$ hold. For each uniformly
overtaking  optimal   pair
    $(x^0,u^0)\in C(\mm{T},\mm{X})\times\fr{U}$
satisfying $\bf{(\psi A)}$, for each unbounded increasing sequence
$(\tau_n)_{n\in\mm{N}}\in \mm{T}^\mm{N}$
 there exists
 $(x^0,u^0,\lambda^0,\psi^0)\in\fr{Z}$ such that
      \begin{equation}
   \label{partlim_2}
       \liminf_{n\to\infty}
       ||\psi^0(\tau_n)A_*(\tau_n)||_\mm{X}=0.
   \end{equation}
 hold.
\end{proposition}
\end{subequations}

  The only differences between the proof of this Proposition and
  Proposition~\ref{3} are the facts that the references to~{($\psi$)}
   are replaced with references to~{($\psi A$)} and the factors
  $A^\eta_{\xi},A_*$ are added to the inequalities of the last strings.

\medskip

     Similarly, we can formulate an analogue to Corollary~\ref{s3} for
this condition: if it is possible to choose matrix maps such that the
product $\psi A^{u}_{z(0)}$ is the solution of an equation
\begin{equation}
 \label{1460}
 \frac{dp}{dt}= b(t,p(t),\psi(t),x(t),\lambda,u(t))
 \end{equation}
for  each $u\in\fr{U}$ for each solution $z=(\psi,x,\lambda)$ of
system~\rref{sys_psi},\rref{sys_x},\rref{sys_lam} with initial
conditions $z(0)\in E_q$, then the corresponding stability of this
component $p$  in system
    \rref{1460},\rref{sys_psii}--\rref{sys_lam}
     implies the result of Corollary~\ref{s3}
(see \cite{my1}).

 The simplest way to account for the a priori information on
stability or for asymptotic estimates of~$\psi$ and its components
is to take $A_*(t)\rav e^{-\lambda t}1_\mm{L}$, where~$\lambda$ is
greater than or equal to all Lyapunov's exponents of the
variable~$\psi$. In particular, in~\cite[Example~10.2]{norv}, the
use of $A_*(t)\rav e^{-t}1_\mm{L}$ in such condition (in contrast to
the standard condition) selects the single extremal.

\section{Cauchy formula for adjoint variable} \label{sec:6}
\ee In the papers~\cite{kr_as1,kr_as2,kr_asD,kr_as,kr_as3}, Aseev and
Kryazhimskii have proposed and proven the analytic
  expression for the values of the adjoint variables.
  This version of the normal form of the Maximum Principle
holds with the explicitly specified adjoint variable providing a
complete set
  of necessary optimality conditions; moreover, the solution of this form of Maximum Principle is uniquely determined by the optimal
  control.
 This approach
generalizes (see~\cite[Sect.~16]{kr_as}, \cite{av}) a number of
transversality conditions; in particular, it is more general than the
conditions that were obtained for linear systems in~\cite{aucl}.

It turns out that  if the function $A_*$ is fundamental matrix of linearized system
  along the
optimal
 trajectory, then,  condition \rref{partlim_1}
 automatically yields this explicit representation for the adjoint variable.

 Let us simplify Proposition~\ref{3_} for such~$A_*$ to weaken the
requirements of \cite[Theorem~ 12.1]{kr_as},\cite[Theorem~
1]{av},\cite[Theorem~2]{kr_as3}, and their corollaries.

\subsection{The case of dominating discount} \label{sec:61}

Let the pair $(x^0,u^0)\in C_{loc}(\mm{T},\mm{X})\times\fr{U}$
  be  uniformly  overtaking  optimal   for problem~\rref{sys}--\rref{opt}.
   Along with it, let us consider the
solution $A_*$ of the Cauchy problem
 $$\frac{d{A_*}(t)}{dt} =\frac{\partial f (t,x^0(t),u^0(t))}{\partial x} A_*(t),\quad
  A_*(0)=1_\mm{L}.$$

Likewise, for each $\xi\in\mm{X}$, let us denote by~$x_{\xi}$  the
solution of~\rref{sys_x} for the initial condition
$x_{\xi}(0)=\xi\in\mm{X}$; let us also consider~$A_{\xi}$, the
solution of the matrix Cauchy problem
 $$\frac{d{A}_{\xi}(t)}{dt}=\frac{\partial
f(t,x_{\xi}(t),u^0(t))}{\partial x} A_{\xi}(t),\quad
A_{\xi}(0)=1_\mm{L}\qquad \forall \xi\in\mm{X}.$$ For each
$T\in\mm{T}$, let us consider
 $$  I_{\xi}(T)\rav\int_0^T
   \frac{\partial g(t,x_{\xi}(t),u^0(t))}{\partial x}\, A_\xi(t)
  \,dt.$$

\begin{proposition}
\label{aff}
 Assume conditions $\bf{(u),(fg),(\partial)}$.
 Let the pair $(x^0,u^0)\in C_{loc}(\mm{T},\mm{X})\times\fr{U}$
  be  uniformly  overtaking  optimal   for problem~\rref{sys}--\rref{opt}.
Let
the map~$I_0$ be bounded and let
$$\displaystyle \lim_{\xi\to 0}
||I_{\xi}-I_0||_C=0.$$ Let $I_*\in\mm{X}$ be a partial limit (the
limit of a subsequence) of
 $I_0(\tau)$ as $\tau\to \infty$.

 Then, there exists a solution
 $(x^0,u^0,\lambda^0,\psi^0)\in\fr{Z}$ of all relations of the Maximum
 Principle \rref{sys_x}--\rref{dob} satisfying
 the transversality condition~\rref{partlim_1}.
 Moreover, $\displaystyle \lambda^0\rav \frac{1}{1+||I_*||_\mm{X}}>0,$
and $\psi^0\in C_{loc}(\mm{T},\mm{X})$ defined by the following rule:
  \begin{equation}
   \label{klass_}
 \psi^0(T)\rav
 \lambda^0 \Big(I_*-\int_0^T
 \frac{\partial g(t,x^0(t),u^0(t))}{\partial x}\,A_*(t)\,dt\Big)
 A_*^{-1}(T) \qquad\forall T\in\mm{T}.
\end{equation}
\end{proposition}

\doc
For each $u\in{\fr{U}}$ and $z\rav{{\ph}[u]}$, let us introduce a
matrix function $A^u$ that is the solution of the matrix equation
\begin{equation}
   \label{Aeta}
 \dot{A}^u(t)=\frac{\partial
f(t,z(t),u(t))}{\partial x}\,A^u(t),\qquad A^u(0)=1_\mm{L}.
\end{equation}
Now, for each solution $(z,u,\lambda^u,\psi^u)\in{\fr{Y}}$
from~\rref{sys_psi_} it follows that
\begin{equation}
   \label{AetaPsi_1}
\frac{d}{dt}\big(\psi^u A^u\big)(t)=-\lambda^u  \frac{\partial
g(t,z(t),u(t))}{\partial x} A^u(t).
\end{equation}

For
\begin{eqnarray}
  E_p\rav\mm{X},\ &\ &  b(t,p,(q_1,q_2,q_3),u)\rav-q_3 \frac{\partial
 g(t,q_2,u)}{\partial x}  q_1, \nonumber \\[-1.5ex]
\label{eq2.14}\\[-1.5ex]
E_q\rav\mm{L}\times\mm{X}\times\mm{R},\ &\
&c(t,p,(q_1,q_2,q_3),u)\rav
 \Big( \frac{\partial
 f(t,q_2,u)}{\partial x}q_1,f(t,q_2,u),0\Big)\nonumber
 \end{eqnarray}
the system~\rref{aa} becomes
system~\rref{AetaPsi_1},\rref{Aeta},\rref{sys_x},\rref{sys_lam};
 now, for $u=u^0$,
 \begin{equation}
   \label{2049}
     \dot{p}=-r \frac{\partial
 g(t,z,u^0(t))}{\partial x} B , \dot{B}=\frac{\partial
 f(t,z,u^0(t))}{\partial x} B ,\dot{z}=f(t,z,u^0(t)),
 \dot{r}=0.
\end{equation}
Solving this system, we obtain
 \begin{eqnarray}
   r(t)=r(0),\quad  z(t)=x_{z(0)}(t),\quad
   B(t)=A_{z(0)}(t)B(0),\nonumber \\
     p(t)=p(0)-r(0)I_{z(0)}(t)B(0).
        \label{1943}
\end{eqnarray}
Let
 $G_0\rav\mm{X}\times\mm{L}\times\mm{X}\times[0,1]$.
 We claim  the partial Lyapunov stability of the component~$p$ in $G_0$
 for $B(0)=1_\mm{L},z(0)=0,r(0)\in[0,1], p(0)\in\mm{D}$. Indeed, by
\rref{1943} it remains to verify that~$I_{\xi}$ is continuous at the
point $\xi=0$, and $I_0$ is bounded; both hold by assumptions.
Therefore, by Proposition~\ref{zvez}, for the certain weight $w^0\in
(Null)(u^0)$, the component $p$ is stable for
 $\fr{L}_{w^0}-$small perturbations of the control~$u^0$.
 This provides condition~{($\psi A$)} for~$A^u$ that were
 defined as we did.

By condition,~$I_*$ is a partial limit; hence, there exists an
unbounded increasing sequence $(\tau_n)_{n\in\mm{N}}\in\mm{T}^\mm{N}$
with property~$I_0(\tau_n)\to
 I_*$. Now, by Proposition~\ref{3_}, there exists
 $(x^0,u^0,\lambda^0,\psi^0)\in\fr{Z}$ with properties~\rref{partlim_2}
 and \rref{partlim_1}.

Substituting $z(0)=0$, $r(0)=\lambda^0,$ $A_0(0)=1_\mm{L}$, and
  $p(0)=\psi^0(0)$ into~\rref{1943} yields
\begin{equation}
   \label{AetaPsi_2}
   \big(\psi^0A_*\big)(T)=\big(\psi^0A_0\big)(T)=p(T)=\psi^0(0)-\lambda^0I_0(T)\qquad\forall T\in\mm{T}.
\end{equation}
Now, substituting~$T=\tau_n$ and passing to the lower limit, from~\rref{partlim_1}, we obtain
$0=\psi^0(0)-\lambda^0I_*$; therefore, from~\rref{AetaPsi_2}
and~\rref{dob} respectively, we have
 $$\displaystyle \psi^0(T)A_*(T)=\lambda^0\big(I_*-I_0(T)\big),\quad
  \lambda^0=\frac{1}{1+||I_*||_\mm{X}}> 0.$$
Using the inverse matrix for~$A_*(T)$,
 we obtain~\rref{klass_}.
 \bo

Let us note that if~$I_*$ is independent of choice of the subsequence
$(\tau_n)_{n\in\mm{N}}$, we automatically obtain the stronger
transversality condition
   \begin{equation}
   \label{lim}
   \lim_{t\to\infty} \psi(t)A_*(t)=0.
\end{equation}
Moreover, since for different $(x^0,u^0,\lambda)$, solutions
of~\rref{AetaPsi_1} differ by a constant, for all
  $(x^0,u^0,\lambda,\psi)\in\fr{Y}$, the products $\psi A_0$
  tend to a finite limit as $t\to\infty$. If \rref{klass_} holds, then
  this limit is equal to zero.
  Hence, to every $(x^0,u^0,\lambda)$ there corresponds at most one~$\psi^0$,
   for which relations
  \rref{sys_x}--\rref{maxH},
  \rref{lim} hold;
now, from~\rref{dob} and \rref{1943} we can reconstruct~$\lambda^0$
uniquely. Thus there exists the unique solution
$(x^0,u^0,\lambda,\psi)\in\fr{Z}$ that satisfies
condition~\rref{lim}, and the following theorem is proved.

\begin{theorem}
\label{aff_}
 Assume conditions $\bf{(u),(fg),(\partial)}$ hold. Let the pair $(x^0,u^0)\in C_{loc}(\mm{T},\mm{X})\times\fr{U}$
  be  uniformly  overtaking  optimal   for
  problem~\rref{sys}--\rref{opt},
and let the limit
   \label{omega2} 
$$   \lim_{t\to\infty,\xi\to 0_\mm{X}} I_{\xi}(t)=\int_{\mm{T}}
 \frac{\partial g(t,x^0(t),u^0(t))}{\partial x}  A_*(t) dt\in\mm{R}$$
be well-defined and finite.

 Then, there exists the unique  solution
 $(x^0,u^0,\lambda^0,\psi^0)\in\fr{Z}$ of all relations of the Maximum
 Principle \rref{sys_x}--\rref{dob} satisfying
 the transversality condition~\rref{lim}.
 Moreover, accurate to the positive factor, we can assume
  \begin{equation}
   \label{klass}
   \lambda^0\rav 1,\qquad
 \psi^0(T)\rav
 \int_{[T,\infty\ra}
 \frac{\partial g(t,x^0(t),u^0(t))}{\partial x}  A_*(t) \,dt\, {A_*^{-1}(T)}\qquad\forall T\in\mm{T}.
\end{equation}
\end{theorem}

From conditions of
\cite[Theorem~2]{kr_as3},\cite[Theorem~12.1]{kr_as},\cite[Theorem~1]{kr_asD},\cite[Theorem~
1]{av} it follows that for some $\alpha,\beta\in\mm{R}_{>0}$ and for
all admissible controls~$u$,  all trajectories~$x$, and all
fundamental matrices~$A^u$, the inequality
  \begin{equation}
   \label{exp}
  \Big|\Big|\frac{\partial g(t,x(t),u(t))}{\partial x}\Big|\Big|\, ||A^u(t)||\leq \beta e^{-\alpha
  t}\qquad
   \forall t\in\mm{T}
\end{equation}
  holds. This is stronger than the conditions of Theorem~2.
Informally, the requirements
of~\cite[Theorem~1]{kr_asD},\cite[Theorem~2]{kr_as3},\cite[Theorem~12.1]{kr_as},
\cite[Theorem~1]{av} boil down to the need for uniform exponential
Lyapunov stability of the product~$\psi A$ along all trajectories of
the system~\rref{1}, while Lyapunov stability of the product $\psi^0 A_0$ along
the optimal solution of the initial control problem is sufficient for
Theorem~2. On the other side, the condition \rref{exp} can be
verified by calculating the Lyapunov exponents of the system of the
Maximum Principle,
see~\cite[Sect.~12]{kr_as},\cite[Sect.~3]{kr_as3},\cite[Sect.~5]{av}.

\subsection{The general case} \label{sec:62}

Let us base on the start of the proof of
Proposition~\ref{aff}, and let us use not~Corollary~\ref{s3} but Remark~\ref{1014}.
 \begin{proposition}
 \label{2190}
 Assume conditions $\bf{(u),(fg),(\partial)}$ hold.
 Let the pair $(x^0,u^0)\in C_{loc}(\mm{T},\mm{X})\times\fr{U}$ be   uniformly  overtaking  optimal
  for problem~\rref{sys}--\rref{opt}.

Now, for an unbounded increasing sequence of times
$(\tau_n)_{n\in\mm{N}}\in \mm{T}^\mm{N}$, there exist:
  \begin{description}
\item[1)] its subsequence $(t_n)_{n\in\mm{N}}\in
   \mm{T}^\mm{N}$;
\item[2)] the sequence of initial conditions $(\zeta_n)_{n\in\mm{N}}\in \mm{X}^\mm{N}$ that converges to~$0_\mm{X}$;
\item[3)] the sequence $(\lambda_n)_{n\in\mm{N}}\in [0,1]^\mm{N}$ that converges to some $\lambda^0\in [0,1]$;
\end{description}
such that if $\psi^0\in C(\mm{T},\mm{X})$ is defined for every $t\in\mm{T}$ by the rule
$$  \psi^0(T)=
  \lim_{n\to \infty}
  \lambda_n
 \int_T^{t_n}
 \frac{\partial g(t,x_{\zeta_n}(t),u^0(t))}{\partial
 x}A_{\zeta_n}(t)\,dt\,A_*^{-1}(T),$$
then, the limit would be uniform on every compact, and   $(x^0,u^0,\lambda^0,\psi^0)_{n\in\mm{N}}\in
   {\fr{Z}}$ would satisfy all  relations of the Maximum Principle
   \rref{sys_x}--\rref{dob}.
   Moreover,
     \begin{equation}
   \label{deimos_1}
  \psi^0(T)=
  \lim_{n\to \infty} \lambda^n
 \int_T^{t_n}
 \frac{\partial g(t,x_{\zeta_n}(t),u^0(t))}{\partial
 x}A_{\zeta_n}(t)\,dt\,A^{-1}_{\zeta_n}(T)\quad
 \forall T\in\mm{T}.
\end{equation}
\end{proposition}
\doc
Indeed, consider the control system $(\dot{p},\dot{q})=(b,c)=a$ from~\rref{eq2.14}. It features the set of controls~$\fr{U}$,
 however, as a system of form~\rref{a}, it defines the control system of form~\rref{1650}
  that is controlled by the elements of~$\td{\fr{U}}$. For such a system, fix the weight~$w^0$ from the formulation of Proposition~\ref{dop}.

By Remark~\ref{1014}, for every $(\tau_n)_{n\in\mm{N}}\in \mm{T}^\mm{N}$, there exist its subsequence $(t_n)_{n\in\mm{N}}\in
   \mm{T}^\mm{N}$ and the sequence
   $(x^n,\eta^n,\lambda^n,\psi^n)_{n\in\mm{N}}\in\td{\fr{Y}}$ of solutions of
    system \rref{sys_x}--\rref{sys_psi}, converging to the certain solution $(x^0,\td{u}^0,\lambda^0,\psi^0)\in
   \td{\fr{Z}}$ of all relations of the Maximum Principle.

Now, for every $n\in\mm{N}$, we can find
   $B_n\in C(\mm{T},\mm{L})$ and $p_n\in C(\mm{T},\mm{X})$ such that
  \begin{equation}
   \label{2244}
   a^n\rav(p_n,B_n,x^n,\lambda^n)\in\td{\fr{A}}[\eta^n],\quad
      B_n(0)=1_\mm{L},\quad p_n(0)=\psi^n(0).
\end{equation}
On the other side, differentiating $\psi^n B_n$
   (as in~\rref{AetaPsi_1}), we check that $(\psi^n B_n,B_n,x^n,\lambda^n)\in\td{\fr{A}}[\eta^n]$. Comparing the initial conditions, we see that
 $   p_n\equiv \psi^n B_n.$

For each $n\in\mm{N}$, for each $t\in\mm{T}$, there exists the solution $a_{n,t}\in C(\mm{T},E)$
of \rref{2049} for the initial conditions
    $a_{n,t}(0)=\beg(t,a^n(t)).$ Note that the last components of
$a_{n,t}$ and~$a^n$ are independent of~$t$; thus, they correspond with~$\lambda^n$. Now we can correctly define the components of the map $t\mapsto\beg(t,a^n(t))$ by the rule
 $$   \beg(t,a^n(t))=\big(\nu_{n}(t),\mu_{n}(t),\xi^n(t),\lambda^n\big)\qquad    \forall t\in\mm{T},n\in\mm{N}.$$
Substituting
 these initial conditions
  into~\rref{1943}, by virtue of equalities~\rref{2244} and $a^n(t)=a_{n,t}(t)$
  for all $n\in\mm{N},t\in\mm{T}$, we obtain
$$   \big(\psi^n(t)B_n(t) ,B_n(t),x^n(t)\big)=a_{n,t}(t)=
   \big(\nu_n(t)-\lambda^nI_{\xi^n(t)}\mu_n(t),
     A_{\xi^n(t)}\mu_n(t),x_{\xi^n(t)}\big).$$ 
Specifically,
      $$\psi^n(t)A_{\xi^n(t)}(t)=\psi^n(t)B^n(t)\mu_n^{-1}(t)=
      \nu_n(t)\mu_n^{-1}(t)-\lambda^n I_{\xi^n(t)}(t).$$
Remark~\ref{1014} provides $\psi^n(t_n)=0$; substituting $t=t_n$, we obtain      $$0=\psi^n(t_n)A_{\xi^n(t_n)}(t_n)=\nu_n(t_n)\mu_n^{-1}(t_n)-\lambda^n(t_n)I_{\xi^n(t_n)}(t_n).$$
Substracting one from another yields
  \begin{equation}
   \label{2314}
      \psi^n(t)A_{\xi^n(t)}(t)=\nu_n(t)\mu_n^{-1}(t)-\lambda^n I_{\xi^n(t)}(t)-
      \nu_n(t_n)\mu_n^{-1}(t_n)+\lambda^n I_{\xi^n(t_n)}(t_n).
\end{equation}

By Remark~\ref{1014},  we have
    $a^n(0)\ravref{2244}
    (\psi^n(0),1_\mm{L},0_\mm{X},\lambda^n)\to (\psi^0(0),1_\mm{L},0_\mm{X},\lambda^0),$
       and   $||\fr{L}_{w^0}(\eta^n)||_{C}\to
     0$
as $n\to\infty$; moreover, Proposition~\ref{dop} yields for all $t\in\mm{T}$
    $$\max_{t\in\mm{T}}||a_{t,n}(0)-a^n(0)||_{E}=
    \max_{t\in\mm{T}}||\beg(t,a^n(t))-\beg(0,a^n(0))||_{E}\leq ||\fr{L}_{w^0}(\eta^n)||_{C}\to 0.$$
Hence uniformly on the whole~$\mm{T}$ as $n\to\infty$, it holds that
  \begin{equation}
   \label{2277}
 a_{t,n}(0)=\big(\nu_{n}(t),\mu_{n}(t),\xi^n(t),\lambda^n\big)\to(\psi^0(0),1_\mm{L},0_\mm{X},\lambda^0)
 \qquad\forall t\in T.
\end{equation}
Whence the theorem of continuous dependence on initial conditions yields the uniformity of the limits             $$\lim_{n\to\infty} \nu_n(\tau)\mu_n^{-1}(\tau)=\psi^0(0),\quad
            \lim_{n\to\infty} \lambda^n I_{\xi^n(\tau)}(t) =\lambda^0
      I_{0}(t)
      \qquad\forall t\in K,\tau\in\mm{T}.
      $$
as $n\to\infty$ for each compact $K\in(comp)(\mm{T})$. Putting here $\tau=t,$ $\tau=t_n$, let us consider
the limit of both sides of~\rref{2314} as $n\to\infty$; thus,
\begin{eqnarray*}
      \psi^0(t)A_{0}(t)=\lim_{n\to\infty} \Big(\psi^0(0)-\lambda^0 I_{0}(t)-\psi^0(0)+
      \lambda^n I_{\xi^n(t_n)}(t_n)\Big)=\\
      \lim_{n\to\infty}\Big(-\lambda^n I_{\xi^n(t_n)}(t)
      +\lambda^n I_{\xi^n(t_n)}(t_n)\big)=
       \lim_{n\to\infty}       \lambda^n \big(I_{\xi^n(t_n)}(t_n)-I_{\xi^n(t_n)}(t)\big).
 \end{eqnarray*}
Multiplying on the right by $A^{-1}_{0}(t)=A_*^{-1}(t)$ and $A_{0}(t)A_{\xi^n(t_n)}(t)$
we obtain our proposition  for $\zeta_n\rav \xi^n(t_n)$.
All necessary convergences are provided by uniformity of limits in~\rref{2277}.
\bo


We say an optimal pair $(x^0,u^0)\in
C_{loc}(\mm{T},\mm{X})\times\fr{U}$ for problem
~\rref{sys}--\rref{opt} is {\it abnormal} if every solution
$(x^0,u^0,\lambda^0,\psi^0)\in
   {\fr{Z}}$ of all relations of the Maximum Principle \rref{sys_x}--\rref{dob}
    satisfies $\lambda^0=0.$

\begin{remark}
 \label{2350}
 Assume conditions $\bf{(u),(fg),(\partial)}$.
 Let the pair $(x^0,u^0)\in C_{loc}(\mm{T},\mm{X})\times\fr{U}$
  be  uniformly  overtaking  optimal pair for problem~\rref{sys}--\rref{opt}
  and let this pair be abnormal.
Then,
 $$ \limsup_{\tau\to\infty, \xi\to 0}||I_\xi(\tau)||_E=
 \limsup_{\tau\to\infty, \xi\to 0}
   \bigg|\bigg| \int_0^\tau
 \frac{\partial g(t,x_{\xi}(t),u^0(t))}{\partial
 x}A_{\xi}(t)\,dt\, \bigg|\bigg|_E=\infty.$$
 \end{remark}
Indeed, if it is wrong, then, the right-hand side of
   \rref{deimos_1} equals zero for $T=0$, i.e., $\psi^0(0)=0_\mm{X},$ which contradicts
   the relation \rref{dob} for $\lambda^0=0$.

\subsection{Monotonous case}
 \label{sec:63}

 Consider the case of when both the right-hand side of the dynamics equation and
 the objective function are monotonous. This case frequently arises in economical
 applications while monotonicity simplifies its examination.
 It seems that the first to note the peculiarities of this case
  and to investigate it were Aseev, Kryazhimskii, and Taras'ev
  in their paper~\cite{kr_as_t}. These were followed
  by papers~\cite{weber},\cite{kr_as1},\cite{kr_as2},
  and the most general case was considered in~\cite{kr_as}.

In Euclidean space~$E'$, let us define binary relations~$\succcurlyeq$ and~$\succ$ by the rules
   $$(\alpha\succcurlyeq \beta)\Leftrightarrow (\alpha-\beta\in\mm{T}^{dim E}),\quad
   (\alpha\succ \beta)\Leftrightarrow (\alpha-\beta\in\mm{R}_{>0}^{dim E})\qquad
    \forall \alpha,\beta\in E'.$$
This allows us to use the symbols $\succcurlyeq$ and $\succ$ to compare vectors and matrices, and vector and matrix functions. For the latter two, $\succcurlyeq$ and $\succ$ allow us to discuss their monotonicity.

\begin{proposition}
\label{mon}
 Assume conditions $\bf{(u),(fg),(\partial)}$ hold.
 Let the pair $(x^0,u^0)\in C_{loc}(\mm{T},\mm{X})\times\fr{U}$
  be  uniformly  overtaking  optimal for problem~\rref{sys}--\rref{opt}.
  Assume for all $x\in\mm{X}$ and for a.a. $t\in\mm{T}$
   there exists a number $d(t,x)\in\mm{R}$  such that the following relation holds:
 $$ \frac{\partial g(t,x,u^0(t))}{\partial x}\succcurlyeq 0_\mm{L},\quad
  \frac{\partial f(t,x,u^0(t))}{\partial x}\succcurlyeq d(t,x)1_\mm{L}.$$
Then, there exists a solution
 $(x^0,u^0,\lambda^0,\psi^0)\in\fr{Z}$ of all relations of the Maximum
 Principle \rref{sys_x}--\rref{dob} satisfying \rref{deimos_1}, and $\psi^0\succcurlyeq 0_\mm{X}.$

If at the same time the pair $(x^0,u^0)$ is normal, then
\begin{equation}
 \label{2455}
   \lambda^0 \limsup_{t\to\infty,\xi\to 0} I_\xi(t)
    \succcurlyeq \psi^0(0) \succcurlyeq \lambda^0
   \lim_{t\to\infty} I_0(t) \succcurlyeq 0_\mm{X}
\end{equation}
 hold, and all limits in \rref{2455} well-defined and finite.
\end{proposition}
\begin{corollary}
\label{mon2}
 Assume conditions $\bf{(u),(fg),(\partial)}$ hold.
 Let the pair $(x^0,u^0)\in C_{loc}(\mm{T},\mm{X})\times\fr{U}$
  be  uniformly  overtaking  optimal for problem~\rref{sys}--\rref{opt}, and
  let this pair be normal.
  Assume for all $x\in\mm{X}$ and for a.a. $t\in\mm{T}$
   there exists a number $d(t,x)\in\mm{R}$  such that the following relation holds:
$$  \frac{\partial g(t,x,u^0(t))}{\partial x}\succ 0_\mm{L},\quad
  \frac{\partial f(t,x,u^0(t))}{\partial x}\succ d(t,x)1_\mm{L}.$$
Then, there exists a solution
 $(x^0,u^0,\lambda^0,\psi^0)\in\fr{Z}$ of all relations of the Maximum
 Principle \rref{sys_x}--\rref{dob} satisfying \rref{deimos_1},\rref{2455}, and
 $
   \psi^0\succ 0_\mm{X}.$
\end{corollary}
\doc
Below, in the proof of Proposition~\ref{mon}, we understand the symbol $\vartriangleright$ as $\succcurlyeq$, and in the proof of Corollary~\ref{mon2}, we understand it as $\succ$.

Fix arbitrary $\xi\in\mm{X},T\in\mm{R}_{>0},\tau\in\la T,\infty\ra$;
 let us show that $A_\xi(\tau)A^{-1}_\xi(T)\vartriangleright 0_\mm{L}$.
  Denote by $F_\xi(t)$ the matrix
   $\frac{\partial f(t,x_\xi(t),u^0(t))}{\partial x}$ for all $t\in[T,\tau]$.
  The diagonal of the map~$F_\xi$ is dominated by a
  function~$M\rav M_{[T,\tau]}^{F_\xi}\in\ct{L}^1_{loc}(\mm{T},\mm{T});$
  then, by condition, $F_\xi+m(t)1_\mm{L}|_{\la T,\tau]}
   \vartriangleright 0_\mm{L}.$ Now, let us consider a solution~$P(t)$ of the equation
      $$\dot{P}=(F_\xi(t)+M(t)1_\mm{L}) P, P(T)=1_\mm{L}, t\geq T;$$
for it, it holds that $P(t) \vartriangleright 0_\mm{L}$ for all
$t\in\la T,\tau].$ Since~$A_\xi$ and~$1_\mm{L}$ commute, the solution
$P$~is the product of two solutions of the equations
      $\dot{Q}=F_\xi(t)Q,\ Q(T)=1_\mm{L}$, and
      $\dot{R}=M(t)1_\mm{L} R,\ R(T)=1_\mm{L}.$
Thus, $$P(\tau)=Q(\tau)R(\tau)=Q(\tau) e^{\int_{T}^\tau
M(t)dt}1_\mm{L}=
      A_\xi(\tau)A^{-1}_\xi(T){e^{\int_{T}^\tau M(t)dt}},$$
and $P(\tau) \vartriangleright 0_\mm{L}$ implies
$A_\xi(\tau)A^{-1}_\xi(T) \vartriangleright 0_\mm{L}$ for all
$\tau\in\la T,\infty\ra$. Now, by monotonicity of matrix product, we
obtain
\begin{equation}
   \label{mmm}
   \frac{dI_\xi(t)}{dt}A_\xi^{-1}(T)=
   \frac{\partial g(t,x_\xi(t),u^0(t))}{\partial
    x}A_\xi(t)A_\xi^{-1}(T) \vartriangleright 0_\mm{X}\quad \forall
    t\in\la T,\infty\ra
\end{equation}
 for all $\xi\in\mm{X},T\in\mm{T};$
specifically, for $T=0$ we have $\frac{dI_\xi(t)}{dt} \vartriangleright
0_\mm{X},$ hence  the functions
    $I_\xi,I_\xi A_\xi^{-1}(T)$ are monotonically increasing
 for all
 $\xi\in\mm{X},T\in\mm{T}.$

By Proposition~\ref{2190}, there exists the solution $(x^0,u^0,\psi^0,\lambda^0)$ of relations of the
Maximum Principle satisfying of formula~\rref{deimos_1} for  certain
sequences $\lambda^n$ and $\xi_n$. However, the expression into the
 limit of~\rref{deimos_1}
    lies in $\mm{L}_{\succcurlyeq 0}$ by~\rref{mmm}. Passing to the limit as $n\to\infty$,
     we obtain $\psi^0\succcurlyeq 0_\mm{X}.$

Suppose the pair $(x^0,u^0)$ is normal; then  $\lambda^0>0$. Since
the function
    $I_\xi$ is monotonically increasing and, by Remark~\ref{2350},
    uniformly bounded in the certain neighborhood $0_\mm{X}$,
     the Lebesgue theorem yields the existence of the finite limits in~\rref{2455}.
   Hence,
%
 $$ {\lambda_0}\limsup_{t\to \infty,\xi\to 0} I_{\xi}(t)\succcurlyeq
  {\lambda_0}\lim_{n\to \infty} I_{\zeta_n}(t_n)
 \ravref{deimos_1}{\psi^0(0)}.
  $$
On the other side, monotonicity of~$I_\xi A_\xi^{-1}(T)$ 
yields
\begin{eqnarray*}
  \frac{1}{\lambda^0}\psi^0(T)\ravref{deimos_1}\lim_{n\to \infty}
  \big(I_{\zeta_n}(t_n)-I_{\zeta_n}(T)\big) A_{\zeta_n}^{-1}(T)\sref{mmm}
  \lim_{n\to \infty} \big(I_{\zeta_n}(t)-I_{\zeta_n}(T)\big)
  A_{\zeta_n}^{-1}(T)\\
  =\big(I_0(t)-I_0(T)\big) A_0^{-1}(T)
  \sref{mmm} 0_\mm{X}\qquad \forall T\in\mm{T},t\in\la T,\infty\ra,
\end{eqnarray*}
 i.e.  $\psi^0\vartriangleright  0_\mm{X}.$
 Moreover, substituting $T=0$ and passing to the limit as $t\to\infty,$ we obtain
 the lower estimate from \rref{2455}.
     \bo

Note that in \cite[Theorem 1]{kr_as_t}, \cite[Theorem 10.1]{kr_as}
the estimate $\psi\succcurlyeq 0_\mm{X}$ ($\psi\succ 0_\mm{X}$) is
made for autonomous systems under less general assumptions;
in~\cite[Theorem 10.1]{kr_as}, the lower estimate from \rref{2455}
was made too (see \cite[(10.17)]{kr_as}). However, in these papers,
the condition $\lambda>0$ was not assumed but proved;
 namely, with the aid of the normal-form stationarity condition,
  the boundedness of  integrals of~\rref{deimos_1} was proved,
  which guaranteed the control was normal.

Let us also note that formula~\rref{klass} was also proved for
biaffine control system for monotonic $\frac{\partial g}{\partial
x}$ (\cite[Theorem~1]{kr_asD},\cite[Theorem~11.1]{kr_as}). It seems,
this result is not a direct consequence of Theorem~\ref{aff_}, which
was proved in that paper.

\medskip

\section{Addendum}
\label{sec:7} \ee
 In the paper, the left endpoint is considered to be
fixed. It seems this condition may be easily discarded, since to do
it,
 it is sufficient to equip the finite horizon optimization problems
 from the proof of Theorem~\ref{1} with the same condition for the left
  endpoint and to provide the boundedness of $x(0)$.

\subsection{Case of $\sigma$-compact-valued map $U$}
The condition~${\bf (u)}$ implies that at every time $t\in\mm{T}$, the controls are chosen from the compact~$U(t)$. Let us weaken this assumption to the following:

%
%

{\it \bf Condition}~${\bf{(u_\sigma)}}:$
     $U$ is a $\sigma$-compact-valued map such that 
  $Gr\,  U  \in       \ct{B}(\mm{T}\times\mm{U})$. 

We shall still assume the conditions ${\bf{(a),(fg)}}$ to hold, and
we shall not change the definition of~${\fr{U}}$. Then, we can assume there exists a nondecreasing sequence $(U^{(r)})_{r\in\mm{N}}$ of  integrally bounded (on each compact subset of~$\mm{T}$) compact-valued maps such that
 $U\equiv\cup_{r\in\mm{N}}U^{(r)}$. Let us assume the uniformly  overtaking optimal control~$u^0$ exists. Then, we may safely assume that $Gr u^0\subset Gr U^{(1)}$.

Repeating the reasonings of~Sect~\ref{sec:11}, for each $r\in\mm{N}$, we can construct sets $\fr{U}^{(r)},\td{\fr{U}}^{(r)}$ and their images for the restriction:
  ${\fr{U}}^{(r)}_n\rav\pi_n(\fr{U}^{(r)}),
  \td{\fr{U}}^{(r)}_n\rav\td{\pi}_n(\td{\fr{U}}^{(r)}).$
  Let us introduce the set  $\td{\fr{U}}$ of all maps~$\eta$ from~$\mm{T}$ into
 the set of Radon probability measures over~$ \mm{U}$ such that
 for every $n\in\mm{N}$ there exists $r=r(\eta,n)\in\mm{N}$ such that
 $\pi_n(\eta)=\eta|_{[0,n]}\in \td{\fr{U}}^{(r)}_n$. The topology of this set is of no use to us, thus we assume it is indiscrete. Note that under our definition, $\td{\delta}(\fr{U})\not\subset \td{\fr{U}}$,
but $u^0\in \td{\delta}(\fr{U}^{(r)})\subset \td{\delta}(\fr{U})$ for all $r\in\mm{N}$

Note that, for all $\eta\in\td{\fr{U}}$, the set $\td{\fr{A}}[\eta]$
is still compact. Indeed, for each $n\in\mm{N}$, it holds that
$\td{\fr{A}}[\eta]|_{[0,n]}\subset\td{\fr{A}}[\td{\fr{U}}^{(r(\eta,n))}]|_{[0,n]}\in(comp)(C([0,n],E));$
all that remains is to use the definition of compact-open topology.
For each $r\in\mm{N}$, denote by $Z^{(r)}$ the pairs
$(\psi,\lambda)\in C(\mm{T},\mm{X})$ such that
$(x^0,u^0,\psi^r,\lambda^r)$ satisfies relations
\rref{sys_x}--\rref{sys_psi},\rref{dob}, and for a.a. $t\in\mm{T}$
instead of~\rref{maxH}, there holds the weaker relation
\begin{equation}\label{maxHr}
\sup_{p\in U^{(r)}(t)}\ct{H}\big(x(t),t,p,\lambda,\psi(t)\big)=
        \ct{H}\big(x(t),t,u^0(t),\lambda,\psi(t)\big).
   \end{equation}
Note that this set is compact (it follows from compactness of $\td{\fr{A}}[\eta]$).
By Theorem~\ref{1}, $Z^{(r)}$ is not empty for all $r\in\mm{N}$. It is easily seen that $Z^{(r')}\subset Z^{(r'')}$ for any
 $r',r''\in\mm{N},(r'\subset r'')$. Then, there exists $(\psi^0,\lambda^0)\in\cap_{r\in\mm{N}} Z^{(r)}$.
 Therefore, for it,~\rref{maxHr} holds for all $r\in\mm{N}$; thus,~\rref{maxH} holds too; hence, $(x^0,u^0,\psi^0,\lambda^0)$
 satisfies all relations \rref{sys_x}--\rref{dob}
  of the Maximum Principle.

For each $r\in\mm{N}$, consider the sequences $(x^n_r,\eta^n_r,\lambda^n_r,\psi^n_r)_{n\in\mm{N}},(t_{n,r})_{n\in\mm{N}}$ from Remark~\ref{1014}. Then, for
the sequence $(x^n_n,\eta^n_n,\lambda^n_n,\psi^n_n)_{n\in\mm{N}}$,
by uniformity of estimate~\rref{to_w}, we have pointwise convergence
of~$\eta_n$ to~$\td{u}_0$; moreover, for each $k\in\mm{N}$ in the interval $[0,k]$
 for this sequence, the convergences from Remark~\ref{1014} hold
   (it is sufficient to consider there topologies with respect to 
   $C([0,k],\mm{X})$,$\td{\fr{U}}^{(k)}$), specifically, it would hold that
    \begin{description}
     \item[1)]
     $(x^n,\eta^n,\lambda^n_n,\psi^n_n)\to (x^0,\td{u}^0,\lambda^0,\psi^0)_{n\in\mm{N}}\in
   C_{loc}(\mm{T},\mm{X})\times \td{\fr{U}}\times \mm{T}\times C_{loc}(\mm{T},\mm{X})$;
    \item[2)] $||\fr{L}_w(\eta^n)||_C\to 0$;
    \item[3)] $\td{J}_{t_n}(\eta)-{J}_{t_n}(u^0)\to 0$, and $\psi^n_n(t_n)=0$ for each
$n\in\mm{N},$ where $t_n\rav t_{n,n}$.
 \end{description}

The verbatim repetition of the proof of Proposition~\ref{3} yields
\begin{proposition}
\label{3in}
 Assume conditions ${\bf (u_\sigma),(fg)}$. For each uniformly   overtaking optimal
pair
    $(x^0,u^0)\in C_{loc}(\mm{T},\mm{X})\times\fr{U}$
satisfying $\bf{(\psi)}$ for each unbounded increasing sequence
$(\tau_n)_{n\in\mm{N}}\in \mm{T}^\mm{N}$
 there exists
 $(x^0,u^0,\lambda^0,\psi^0)$ such that
{ the relations of the Maximum Principle \rref{sys_x}--\rref{dob},
and
 the transversality condition~\rref{partlim} hold.}
\end{proposition}

Since our case is more general, starting with Sect.~\ref{sec:3}, the
references to $\bf(u)$ ought to be replaced with $\bf(u_\sigma)$, and
the results of Sect.~\ref{sec:2} ought to be replaced with their
respectful analogues.

\subsection{On uniformly sporadically catching up controls.}
\begin{definition} We say that a control $u^0\in\fr{U}$ is
   uniformly sporadically catching up optimal  if
   for  every  $\epsi,T \in\mm{R}_{>0}$
    there
   exists  $t\in[T,\infty\ra$  such  that
 $
   J_t(u^0)\geq J_t(u)-\epsi $
 holds for all $u\in\fr{U}$.
\end{definition}
Note that for each uniformly sporadically catching up
   optimal control,     there
   exists an unbounded monotonically increasing sequence
    $(\tau_n)_{n\in\mm{N}}\in\mm{T}^\mm{N}$ and
a    function $\omega^0\in\Omega$
    such  that
$$ J_{\tau_n}(u^0)\geq J_{\tau_n}(u)-\omega^0(\tau_k)\qquad\forall u\in
 \fr{U},k,n\in\mm{N},k<n.$$
We call such control a $\tau$-sporadically catching up optimal.

Now, if we consider the sequence $(\tau_n)_{n\in\mm{N}}$ everywhere
defined and understand the optimality in the above sense, then all
statements, starting with Theorem~\ref{1}, hold.
In particular, we can rewrite Proposition~\ref{aff} and Theorem~\ref{aff_}
 in the following way:
\begin{theorem}
 Assume conditions $\bf{(u_\sigma),(fg),(\partial)}$ hold. Let the
pair $(x^0,u^0)\in C_{loc}(\mm{T},\mm{X})\times\fr{U}$ be $\tau$-sporadically catching up optimal for
problem~\rref{sys}--\rref{opt}. Let $I_*\in\mm{X}$ be the limit  of
 $I_\xi(\tau_n)$ as $n\to \infty,\xi\to 0_\mm{X}.$

  Then, there exists the unique solution
 $(x^0,u^0,\lambda^0,\psi^0)$ of all relations of the Maximum
 Principle \rref{sys_x}--\rref{dob} and
 the transversality condition~\rref{partlim_2}.
 Moreover, accurate to the positive factor, we can assume
 $$
   \lambda^0\rav 1,\quad
 \psi^0(T)\rav
 \Big(I_*-\int_0^T
 \frac{\partial g(t,x^0(t),u^0(t))}{\partial x}\,A_*(t)\,dt\Big)
 A_*^{-1}(T) \qquad\forall T\in\mm{T}.$$
\end{theorem}
\begin{acknowledgements}
I would like to express my gratitude to  A.~G.~Chentsov,
A.~M.~Tarasyev, N.~Yu.~Lukoyanov, and Yu.~V.~Averboukh for valuable
discussion in course of writing this article. Special thanks to
Ya.~Salii for the translation.
\end{acknowledgements}

\bibliographystyle{spmpsci}      


\end{document}